\documentclass[onefignum,onetabnum]{siamart190516}

\usepackage[utf8]{inputenc} 
\usepackage{amsmath,amssymb}
\allowdisplaybreaks

\usepackage{algorithm}
\usepackage{algorithmicx}
\usepackage{algpseudocode}

\usepackage{graphicx}
\usepackage{subcaption}
\usepackage{afterpage}
\usepackage{pdflscape}

\usepackage{booktabs}
\usepackage{array}
\newcommand{\bmcol}[2]{\multicolumn{1}{#1}{\textbf{#2}}}
\newcommand{\mcol}[2]{\multicolumn{1}{#1}{#2}}

\newcommand{\tableemph}{\textbf}

\renewcommand{\vec}{\mathbf}
\newcommand{\T}{\textsf{T}}

\usepackage{xspace}
\newcommand{\HSCG}{\hyperref[alg:hs_pcg]{{HS-CG}}\xspace}
\newcommand{\CGCG}{\hyperref[alg:cg_pcg]{\textsc{CG-CG}}\xspace}
\newcommand{\MCG}{\hyperref[alg:m_pcg]{\textsc{M-CG}}\xspace}
\newcommand{\PRCG}[1][]{\hyperref[alg:pr_pcg]{\textsc{PR-CG}#1}\xspace}
\newcommand{\GVCG}{\hyperref[alg:gv_pcg]{\textsc{GV-CG}}\xspace}
\newcommand{\PipeMCG}[1][]{\hyperref[alg:pipe_m_pcg]{\textsc{pipe-P-M-CG}#1}\xspace}
\newcommand{\PipePCG}[1][]{\hyperref[alg:pipe_p_pcg]{\textsc{pipe-P-CG}#1}\xspace}
\newcommand{\PipePRMCG}{\hyperref[alg:pipe_m_pcg]{\textsc{pipe-M-CG}}\xspace}
\newcommand{\PipePRCG}{\hyperref[alg:pipe_pr_pcg]{\textsc{pipe-PR-CG}}\xspace}

\newcommand{\INIT}{\hyperref[alg:initialize]{\textsc{initialize}()}}

\newcommand{\Cgr}{\mathrm{C}_{\text{gr}}}
\newcommand{\Tmv}{\mathrm{T}_{\text{mv}}}
\newcommand{\Ttmv}{\mathrm{T}_{2\text{mv}}}
\newcommand{\Cmv}{\mathrm{C}_{\text{mv}}}

\usepackage{lipsum}
\usepackage{amsfonts}
\ifpdf
  \DeclareGraphicsExtensions{.eps,.pdf,.png,.jpg}
\else
  \DeclareGraphicsExtensions{.eps}
\fi

\newsiamremark{remark}{Remark}
\newsiamremark{hypothesis}{Hypothesis}
\crefname{hypothesis}{Hypothesis}{Hypotheses}
\newsiamthm{claim}{Claim}

\headers{Predict-and-recompute Conjugate Gradients}{Tyler Chen, Erin C. Carson}

\title{Predict-and-recompute conjugate gradient variants%
\thanks{Submitted to the journal's Methods and Algorithms for Scientific Computing section July 23,2019; accepted for publication (in revised form) June 29, 2020; published electronically October 7,2020.
\funding{%
The first author is supported by the National Science Foundation Graduate Research Fellowship Program under Grant No. DGE-1762114. 
Any opinions, findings, and conclusions or recommendations expressed in this material are those of the author and do not necessarily reflect the views of the National Science Foundation.
The second author is supported by the Charles University Primus program project PRIMUS/19/SCI/11.
This work was also facilitated through the use of advanced computational, storage, and networking infrastructure provided by the Hyak supercomputer system funded by the STF at the University of Washington.
}}
}

\author{Tyler Chen\thanks{Department of Applied Mathematics, University of Washington, Seattle (\href{mailto:chentyl@uw.edu}{\texttt{chentyl@uw.edu}})}%
\and
Erin C. Carson\thanks{Faculty of Mathematics and Physics, Charles University (\href{mailto:carson@karlin.mff.cuni.cz}{\texttt{carson@karlin.mff.cuni.cz}})}%
}

\usepackage{amsopn}

\usepackage{environ}
\RenewEnviron{ecc}{\iffalse\expandafter\BODY\fi}
\RenewEnviron{tc}{\iffalse\expandafter\BODY\fi}

\ifpdf
\hypersetup{
  pdftitle={Predict-and-recompute conjugate gradient variants},
  pdfauthor={Tyler Chen and Erin C. Carson}
}
\fi

\begin{document}

\maketitle

\begin{abstract}
The standard implementation of the conjugate gradient algorithm suffers from communication bottlenecks on parallel architectures, due primarily to the two global reductions required every iteration.
In this paper, we study conjugate gradient variants which decrease the runtime per iteration by overlapping global synchronizations, and in the case of pipelined variants, matrix-vector products.
Through the use of a predict-and-recompute scheme, whereby recursively-updated quantities are first used as a predictor for their true values and then recomputed exactly at a later point in the iteration, these variants are observed to have convergence behavior nearly as good as the standard conjugate gradient implementation on a variety of test problems.
We provide a rounding error analysis which provides insight into this observation.
It is also verified experimentally that the variants studied do indeed reduce the runtime per iteration in practice and that they scale similarly to previously-studied communication-hiding variants.
Finally, because these variants achieve good convergence without the use of any additional input parameters, they have the potential to be used in place of the standard conjugate gradient implementation in a range of applications.
\end{abstract}

\begin{keywords}
Krylov subspace methods,
Conjugate Gradient,
Parallel algorithms,
Numerical algorithms
\end{keywords}

\begin{AMS}
65F10,
65N12,
65Y05
\end{AMS}

\section{Introduction and background}

The conjugate gradient method (CG), due to Hestenes and Stiefel \cite{hestenes_stiefel_52}, is a widely-used method for solving a linear system of equations \( \vec{A} \vec{x} = \vec{b} \), when \( \vec{A} \in \mathbb{R}^{n\times n} \) is a large symmetric positive definite matrix.
While the low storage costs and low number of floating point operations per iteration make CG an attractive choice for solving very large sparse systems, the dependency structure of the standard conjugate gradient algorithm given in~\cite{hestenes_stiefel_52}, which we call \HSCG and present in \cref{alg:hs_pcg}, requires that nearly every computation be done in sequence.
In particular, it requires two inner products and one (typically sparse) matrix-vector product per iteration, none of which can occur simultaneously.
On distributed memory parallel systems, this results in a communication bottleneck \cite{demmel_hoemmen_mohiyuddin_yelick_07,ashby_ghysels_heirman_vanroose_12,ballard_carson_demmel_hoemmen_knight_schwartz_14,dongarra_heroux_luszczek_16}.
This is because inner products require a costly global reduction involving communication between all nodes, and the matrix-vector product, even if sparse or structured, requires some level of communication between individual nodes.

To address this bottleneck, many \emph{mathematically equivalent} variants of the CG method have been introduced; see for instance \cite[etc.]{rosendale_83, saad_85, chronopoulos_87, meurant_87,saad_89,chronopoulos_gear_89,strzodka_goddeke_06,ghysels_vanroose_14,eller_gropp_16,cools_vanroose_17,cornelis_cools_vanroose_19}.
Broadly speaking, these variants rearrange the \HSCG algorithm such that the communication occurs less frequently or is overlapped with other computations.
As a result, the time per iteration of these methods may be reduced on parallel machines in certain settings.
The algorithmic variants of CG presented in this paper are most closely related to the communication-hiding variants which we call \MCG, \CGCG, and \GVCG (pipelined CG), introduced by Meurant \cite{meurant_87}, Chronopoulos and Gear \cite{chronopoulos_gear_89}, and Ghysels and Vanroose \cite{ghysels_vanroose_14}, respectively. These communication-hiding variants maintain the iteration structure of \HSCG but add auxiliary vectors such that the computations within an iteration can be arranged so that expensive ones can occur simultaneously, effectively ``hiding'' the cost of communication. 
A comparison of these variants is given in \cref{table:comparison} and full descriptions are given in \cref{sec:previous}.

\begin{algorithm}[H]
\caption{Hestenes and Stiefel Conjugate Gradient (preconditioned)}\label{alg:hs_pcg}
\fontsize{10}{10}\selectfont
\begin{algorithmic}[1]
\Procedure{\HSCG}{$\vec{A}$, $\vec{M}$, $\vec{b}$, $\vec{x}_0$}
\State \INIT
\For {\( k=1,2,\ldots \)}
    \State \( \vec{x}_k = \vec{x}_{k-1} + \alpha_{k-1} \vec{p}_{k-1} \)
    \State \( \vec{r}_k = \vec{r}_{k-1} - \alpha_{k-1} \vec{s}_{k-1} \)%
    ,~ \( \tilde{\vec{r}}_k = \vec{M}^{-1} \vec{r}_k \)
    \State \( \nu_{k} = \langle \tilde{\vec{r}}_k, \vec{r}_k \rangle \)
    \State \( \beta_k = \nu_k / \nu_{k-1} \)
    \State \( \vec{p}_k = \tilde{\vec{r}}_k + \beta_k \vec{p}_{k-1} \)
    \State \( \vec{s}_k = \vec{A} \vec{p}_k \)
    \State \( \mu_k = \langle \vec{p}_k,\vec{s}_k \rangle \)
    \State \( \alpha_k = \nu_k / \mu_k \)
   \EndFor
\EndProcedure
\end{algorithmic}
\end{algorithm}

Notably, however, many communication-hiding variants suffer from numerical problems due to this rearranging of computations within each iteration.
Indeed, it is well known that CG is particularly sensitive to rounding errors and any modification to the \HSCG algorithm can have a significant effect on numerical behavior; see \cite{greenbaum_97,meurant_strakos_06} for summaries of CG and Lanczos in finite precision.
Specifically, both the rate of convergence (the number of iterations to reach a given level of accuracy) and the maximal attainable accuracy of any algorithmic variant of CG may be \emph{severely} impacted by carrying out computations in finite precision.

These effects are particularly pronounced in pipelined methods, such as \GVCG, because the additional auxiliary recurrences can cause rounding errors to be amplified \cite{carson_rozloznik_strakos_tichy_tuma_18}.
For example, as shown in \cref{sec:numerical}, there are many problems for which the final accuracy attainable by \GVCG is orders of magnitude worse than that attainable by \HSCG.
As a result, the practical use of some of the pipelined variants is potentially limited because the algorithms may fail to reach an acceptable level of accuracy, or require so many iterations to do so that there is no benefit to the overall runtime.
While there have been many approaches to improving the numerical properties of these variants, such as residual replacement \cite{ghysels_vanroose_14,cools_yetkin_agullo_giraud_vanroose_18} and the use of shifts in auxiliary recurrences \cite{cools_vanroose_17}, these strategies typically require certain parameters to be selected ahead of time based on user intuition or rely on heuristics for when and how to apply corrections.
Moreover, as in the case of residual replacement, these strategies sometimes result in delayed convergence.

In this paper, we present a communication-hiding variant similar to \MCG which requires a single global synchronization per iteration. Like the \MCG algorithm, our variants employ a predict-and-recompute scheme which recursively computes quantities as a \emph{predictor} for their true values and then \emph{recomputes} them later in the iteration. 
We then introduce ``pipelined'' versions of both this variant and of \MCG which allow the computation and communication from the matrix-vector product and preconditioning steps to be overlapped with the global reduction associated with the inner products. 
We demonstrate numerically that the convergence behavior of our pipelined variants is comparable to that of \HSCG despite overlapping all matrix products with global reductions as in  other pipelined variants such as \GVCG.
This observation is complemented by a rounding error analysis, which provides new insight into the way the predict-and-recompute schemes can improve the numerical behavior. 
More broadly, this work demonstrates that it is reasonable to add extra computation to improve the numerical properties of a communication-hiding CG variant, provided that that the computation is done in a way which does not affect the communication pattern of the algorithm.
All of the algorithms introduced in this paper require exactly the same inputs as \HSCG and therefore require no additional tuning by the end user.

Our primary contributions are as follows. 
First, we provide rounding error analyses of the variants introduced in this paper.
In particular, we derive expressions for the residual gap, which is the quantity that dictates the maximum attainable accuracy, as well as for the error in the three-term Lanczos recurrence, which gives some inidcation of the rate of convergence.
Our analysis provides a rigorous explanation to a conjecture made by Meurant in \cite{meurant_87} that the predict-and-recompute scheme resolves a potential source of instability.
In addition, we provide a range of numerical experiments to further support the claim that our pipelined variants can significantly improve the rate of convergence and ultimately attainable accuracy versus existing pipelined variants.
In particular, we note that in the numerical experiments in \cref{sec:numerical}, our variants appear to converge similarly \HSCG without any additional input parameters.
As such, they have the potential to be used as black box solvers wherever \HSCG is used.
Finally, we demonstrate through a strong scaling experiment that the new variants can reduce the time per iteration in a parallel setting. 

Unless otherwise stated, matrices should be assumed to be of size \( n\times n \) and vectors of size \( n\times 1 \).
The transpose of a matrix is denoted with the superscript \( \T \), and the inverse of the transpose denoted with the superscript \( -\T \).
The standard Euclidean inner product and corresponding spectral/operator norm are respectively denoted \( \langle \cdot,\cdot \rangle \) and \( \|\cdot\| \), and the inner product and norm induced by a positive definite matrix \( \vec{B} \) are denoted \( \langle \cdot, \cdot \rangle_\vec{B} \) and \( \|\cdot\|_\vec{B} \).
While much of the theory about the conjugate gradient algorithm applies to complex systems, we consider real systems for convenience.

The variants studied in this paper are all compatible with a preconditioner \( \vec{M}^{-1} = \vec{R}^{-\T}\vec{R}^{-1} \), which allows the algorithms to implicitly solve the system \( \vec{R}^{-\T} \vec{A} \vec{R}^{-1} \vec{y} = \vec{R}^{-\T} \vec{b} \) and then set \( \vec{x} = \vec{R}^{-1} \vec{y} \), using only products with \( \vec{M}^{-1} \); i.e. the \( \vec{R}^{-1} \) factors \emph{need not be known}.
As summarized in \cref{table:comparison}, the preconditioned variants require some additional memory to store the preconditioned residual and any related auxiliary vectors, and updating these additional vectors has the potential to introduce additional rounding errors.
%
Throughout this paper, we use a tilde (`` \(\sim\) '') above a vector to indicate that, in exact arithmetic, the tilde vector is equal to the preconditioner applied to the non-tilde vector; i.e., \( \tilde{\vec{r}}_k = \vec{M}^{-1}\vec{r}_k \), \( \tilde{\vec{s}}_k = \vec{M}^{-1}\vec{s}_k \), etc.

\newcommand{\costs}[7]{#5 (+#6) & #1 (+#2) & #3 & #7}

\begin{table}\centering
\begin{tabular}{lcccc} \toprule 
    \bmcol{l}{variant} &
    \bmcol{c}{mem.}&
    \bmcol{c}{vec.}&
    \bmcol{c}{scal.}&
    \bmcol{c}{time} 
    \\ \midrule
    \HSCG      & \costs{3}{0}{2}{1}{4}{1}{\(  2 \, \Cgr + \Tmv + \Cmv  \)} \\
    \midrule 
    \CGCG      & \costs{4}{0}{2}{1}{5}{1}{\( \Cgr + \Tmv + \Cmv \)} \\
    \MCG       & \costs{3}{1}{3}{1}{4}{2}{\( \Cgr + \Tmv + \Cmv \)} \\
    \PRCG      & \costs{3}{1}{4}{1}{4}{2}{\( \Cgr + \Tmv + \Cmv \)} \\
    \midrule 
    \GVCG      & \costs{6}{2}{2}{1}{7}{3}{\( \max( \Cgr, \Tmv + \Cmv ) \)} \\
    \PipePRMCG & \costs{5}{3}{3}{2}{6}{4}{\( \max( \Cgr, \Ttmv + \Cmv ) \)} \\
    \PipePRCG  & \costs{5}{3}{4}{2}{6}{4}{\( \max( \Cgr, \Ttmv + \Cmv ) \)} \\
    \bottomrule
\end{tabular}
    \caption{
        Summary of costs for various conjugate gradient variants. Values in parenthesis are the additional costs for the preconditioned variants.
        The quantity \textbf{mem.} gives the number of vectors stored, 
        \textbf{vec.} gives the number of vector updates (\texttt{AXPY}s) per iteration, 
        \textbf{scal.} gives the number of inner products per iteration, and 
        \textbf{time} gives the dominant costs (ignoring vector updates and inner products). 
        \( \Cgr \) is the time spent on communication for a global reduction. \( \Tmv \) and \( \Cmv \) are the times spent on computation and communication, respectively, for a matrix-vector product which depend on the method of matrix multiplication (for instance a dense matrix has \( \Cmv = \Cgr \)). 
        \( \Ttmv \) is the cost of computing two matrix vector products simultaneously, which is less than \( 2\,\Tmv \) if implemented in an efficient way.
        Note that in this abstraction we assume that the communication time is independent of the size of messages sent.  
    }
    \label{table:comparison}
\end{table}

\section{Derivation of new variants}

In this sec\-tion we de\-scribe \PRCG, a new communication-hiding variant which requires only one global synchronization point per iteration. 
This variant is similar to \MCG, introduced by Meurant in \cite{meurant_87}, and the relationship between the two algorithms is discussed.

Then, in the same way that \GVCG is obtained from \CGCG, we ``pipeline'' \PRCG to overlap the matrix-vector product with the inner products.
The order in which operations are done in the pipelined version of \PRCG allows for a vector quantity to be recomputed using an additional matrix-vector product, giving the pipelined predict-and-recompute variant \PipePRCG.
Since this matrix-vector product can occur at the same time as the other matrix-vector product and as the inner products, the communication costs per iteration are not increased. 
For comparison, we also derive a pipelined version of Meurant's \MCG method, which we call \PipePRMCG. 

\cref{table:comparison} provides a comparison between some commonly used communication avoiding variants and the newly introduced variants.
It should be noted that although the number of matrix-vector products and inner products in \PipePRMCG and \PipePRCG are increased, most of this work can be done locally, and they have the same dominant communication costs as \GVCG.

\subsection{A simple communication-hiding variant}

Like the derivation of \MCG in \cite{meurant_87} and the variants introduced in \cite{johnsson_84,rosendale_83,saad_85,saad_89}, we derive \PRCG by substituting recurrences into the inner product \(  \nu_k = \langle \tilde{\vec{r}}_k, \vec{r}_k \rangle \).
This allows us to obtain an equivalent expression for the inner product involving quantities which are known earlier in the iteration.

To this end, we first recall that given \( \tilde{\vec{s}}_k = \vec{M}^{-1}\vec{s}_k \), we have
\begin{align}
    \tilde{\vec{r}}_k = \tilde{\vec{r}}_{k-1} - \alpha_{k-1} \tilde{\vec{s}}_{k-1}.
\end{align}
Then, by substituting the recurrences for \( \tilde{\vec{r}}_k \) and \( \vec{r}_k \) from \cref{alg:hs_pcg} into \( \nu_k = \langle \tilde{\vec{r}}_k, \vec{r}_k \rangle \), we can write
\begin{align}
    \nu_k 
    &= \langle \tilde{\vec{r}}_{k-1} - \alpha_{k-1} \tilde{\vec{s}}_{k-1}, \vec{r}_{k-1} - \alpha_{k-1} \vec{s}_{k-1} \rangle
    \tag*{}\\&= \langle \tilde{\vec{r}}_{k-1}, \vec{r}_{k-1} \rangle - \alpha_{k-1} \langle \tilde{\vec{r}}_{k-1}, \vec{s}_{k-1} \rangle 
    - \alpha_{k-1} \langle \tilde{\vec{s}}_{k-1}, \vec{r}_{k-1} \rangle + \alpha_{k-1}^2 \langle \tilde{\vec{s}}_{k-1}, \vec{s}_{k-1} \rangle. \label{eqn:nuk}
\end{align}

Since \( \vec{M} \), and therefore \( \vec{M}^{-1} \), are symmetric, \( \langle \tilde{\vec{s}}_{k-1}, \vec{r}_{k-1} \rangle = \langle \vec{s}_{k-1}, \tilde{\vec{r}}_{k-1} \rangle \).
Thus,
\begin{align}
   \nu_k 
    = \nu_{k-1} - 2 \alpha_{k-1} \langle \tilde{\vec{r}}_{k-1}, \vec{s}_{k-1} \rangle + \alpha_{k-1}^2 \langle \tilde{\vec{s}}_{k-1}, \vec{s}_{k-1} \rangle. 
    \label{eqn:nuk_simplified}
\end{align}
Once \( \vec{s}_k = \vec{A} \vec{p}_k \) and \( \tilde{\vec{s}}_k = \vec{M}^{-1} \vec{s}_k \) have been computed, we can simultaneously compute the three inner products,
\begin{align}
    \mu_k = \langle \vec{p}_k,\vec{s}_k \rangle
    ,&&
    \sigma_k = \langle \tilde{\vec{r}}_k,\vec{s}_k \rangle
    ,&&
    \gamma_k = \langle \tilde{\vec{s}}_k,\vec{s}_k \rangle. 
    \label{eqn:inner_products}
\end{align}

A variant using a similar expression for \( \nu_k \) was suggested in \cite{mcmanus_johnson_cross_99} and briefly mentioned in \cite{carson_rozloznik_strakos_tichy_tuma_18}. However, one term of their formula for \( \nu_k \) has a sign difference from \cref{eqn:nuk_simplified}, and no numerical tests or rounding error analysis were provided.

\begin{algorithm}
\caption{Predict-and-recompute conjugate gradient}\label{alg:pr_pcg}
\fontsize{10}{10}\selectfont
\begin{algorithmic}[1]
\Procedure{\PRCG}{$\vec{A}$, $\vec{M}$, $\vec{b}$, $\vec{x}_0$}
\State \INIT
\For {\( k=1,2,\ldots \)}
    \State \( \vec{x}_k = \vec{x}_{k-1} + \alpha_{k-1} \vec{p}_{k-1} \)
    \State \( \vec{r}_k = \vec{r}_{k-1} - \alpha_{k-1} \vec{s}_{k-1} \)%
    ,~ \( \tilde{\vec{r}}_k = \tilde{\vec{r}}_{k-1} - \alpha_{k-1} \tilde{\vec{s}}_{k-1} \)
    \State \( \nu'_{k} = \nu_{k-1} - 2 \alpha_{k-1}\sigma_{k-1} + \alpha_{k-1}^2 \gamma_{k-1}  \) \label{pr_pcg:nukp}
    \State \( \beta_{k} = \nu'_{k} / \nu_{k-1} \)
    \State \( \vec{p}_k = \tilde{\vec{r}}_k + \beta_k \vec{p}_{k-1} \)
    \State \( \vec{s}_k = \vec{A} \vec{p}_{k} \)%
    ,~ \( \tilde{\vec{s}}_k = \vec{M}^{-1} \vec{s}_k \)
    \State \( \mu_k = \langle \vec{p}_k,\vec{s}_k \rangle \)%
    ,~ \( \sigma_k = \langle \tilde{\vec{r}}_k, \vec{s}_k\rangle \)
    ,~ \( \gamma_k = \langle \tilde{\vec{s}}_k, \vec{s}_k \rangle \)%
    ,~ \( \nu_k = \langle \tilde{\vec{r}}_k, \vec{r}_k \rangle \) \label{pr_pcg:inners} 
    \State \( \alpha_k = \nu_k / \mu_k \)
\EndFor
\EndProcedure
\end{algorithmic}
\end{algorithm}

As shown in \cref{fig:no_recompute}, if \cref{eqn:nuk_simplified} is used to recursively update \( \nu_k \), the final accuracy can be severely impacted.
This phenomenon was observed previously in \cite{rosendale_83,johnsson_84,saad_85}, all of which study variants which use expressions for \( \nu_k \) similar to that in \MCG.
This catastrophic loss of accuracy is caused by the updated value of \( \nu_k \) becoming negative.
In \cite{meurant_87}, Meurant suggests using the recursively-updated value of \( \nu_k \) as a \emph{predictor} for the true value in order to update any vectors required for the algorithm to proceed, and then to \emph{recompute} \( \nu_k = \langle \tilde{\vec{r}}_k,\vec{r}_k \rangle \) at the same time as the other inner products.
We observe experimentally that using this strategy effectively brings the ultimately attainable accuracy to a level similar to that of \HSCG. 
Our rounding error analysis in \cref{sec:rounding} explains this improvement. 
This algorithm, denoted \PRCG, is given in \cref{alg:pr_pcg}.
Note that we use a prime (``~\('\)~'') to distinguish the recursively updated quantity \( \nu'_k \) from the explicitly computed quantity \( \nu_k \).

It is not hard to see that \cref{eqn:nuk_simplified} can be simplified further. 
Indeed, observing that \( \langle \tilde{\vec{s}}_k, \vec{r}_k \rangle = \langle \vec{p}_k, \vec{s}_k \rangle \), which follows from the \( \vec{A} \)-orthogonality of \( \vec{p}_k \) and \( \vec{p}_{k-1} \), and using the expression \( \alpha_k = \nu_k / \mu_k = \langle \tilde{\vec{r}}_k, \vec{r}_k \rangle / \langle \vec{p}_k, \vec{s}_k \rangle \), we can write
\begin{align*}
    \nu_k = - \langle \tilde{\vec{r}}_{k-1}, \vec{r}_{k-1} \rangle + \alpha_k^2 \langle \tilde{\vec{s}}_{k-1}, \vec{s}_{k-1} \rangle
     = - \nu_{k-1} + \alpha_k^2 \langle \tilde{\vec{s}}_{k-1}, \vec{s}_{k-1} \rangle.
\end{align*}
This is exactly the expressions studied in \cite{meurant_87} and used by \MCG; i.e., \MCG can be obtained by replacing line~\ref{pr_pcg:nukp} of \PRCG with \( \nu_k' = -\nu_k + \alpha_k^2 \gamma_k \).
A similar expression, \( \nu_k = -\alpha_k\mu_k + \alpha_k^2 \gamma_k \), was studied in \cite{strzodka_goddeke_06}.

The natural question is then: What is the advantage of the expression for \( \nu_k \) used in \PRCG over the expression used in \MCG?
The incomplete response is that experimentally it seems to work a bit better, especially in the pipelined versions of these variants.
Intuitively, this is because the less simplification that occurs, the closer the variant is to \HSCG.
Of course then it is reasonable to wonder if we should use \cref{eqn:nuk_simplified} instead of \cref{eqn:nuk}.
The again incomplete response is that this change doesn't seem to significantly influence the numerical behavior. 
This is perhaps not particularly surprising; for instance, in the case of the unpreconditioned variants, \cref{eqn:nuk} and \cref{eqn:nuk_simplified} are equivalent even in finite precision. 
An extended discussion on this topic is included in \cref{sec:other_variants}; a full theoretical explanation remains future work. 

\subsection{New pipelined variants}

Recall that our goal is to be able to compute the matrix-vector product and inner products simultaneously.
To this end, note that in both \MCG and \PRCG 
we have the recurrence \( \vec{p}_k = \tilde{\vec{r}}_k + \beta_k \vec{p}_{k-1} \). Thus, defining \( \vec{w}_k = \vec{A} \tilde{\vec{r}}_k \), we can write
\begin{align}
    \vec{s}_k = \vec{A} \vec{p}_k
    = \vec{A} \tilde{\vec{r}}_k + \beta_k  \vec{A} \vec{p}_{k-1}
    = \vec{w}_k + \beta_k \vec{s}_{k-1}.
\end{align}
Similarly, defining \(\vec{u}_k = \vec{A} \tilde{\vec{s}}_k \),
\begin{align}
    \vec{w}_k = \vec{A} \tilde{\vec{r}}_k
    &= \vec{A} \tilde{\vec{r}}_{k-1} - \alpha_{k-1} \vec{A} \tilde{\vec{s}}_{k-1}
    = \vec{w}_{k-1} - \alpha_{k-1} \vec{u}_{k-1}.
\end{align}
Using these recurrences allows us to compute the product \( \vec{u}_k = \vec{A} \tilde{\vec{s}}_k \) at the same time as all of the inner products.

To move the preconditioning step, we define \( \tilde{\vec{w}}_k = \vec{M}^{-1} \vec{w}_k \) so that
\begin{align}
    \tilde{\vec{s}}_k = \vec{M}^{-1} \vec{s}_k
    &= \vec{M}^{-1} \vec{w}_k + \beta_k \vec{M}^{-1}\vec{s}_{k-1}
    = \tilde{\vec{w}}_k + \beta_k \tilde{\vec{s}}_{k-1},
\end{align}
and define \( \tilde{\vec{u}}_k = \vec{M}^{-1}\vec{u}_k \) so that
\begin{align}
    \tilde{\vec{w}}_k = \vec{M}^{-1} \vec{w}_k
    &= \vec{M}^{-1} \vec{w}_{k-1} - \alpha_{k-1} \vec{M}^{-1} \vec{u}_{k-1}
    = \tilde{\vec{w}}_{k-1} - \alpha_{k-1} \tilde{\vec{u}}_{k-1}. 
\end{align}

The recurrences above could be used to implement a pipelined variant which, in each iteration, overlaps the matrix product and preconditioning steps with the global reduction.
However, like \GVCG, this variant suffers from delayed convergence and reduced final accuracy compared to \HSCG, \MCG, and \PRCG.
To address this, we observe that \( \vec{w}_k = \vec{A} \tilde{\vec{r}}_k \) and \( \tilde{\vec{w}}_k = \vec{M}^{-1}\vec{w}_k \) can be \emph{recomputed} at the same time as the other matrix-vector product and all inner products are being computed.
Thus, in the same way we use the recursively-updated value of \( \nu_k \) as a predictor for the true value, we can use the recursively-updated value of \( \vec{w}_k \) as a predictor for the true value in order to update other vector quantities, and then update the value of \( \vec{w}_k \) later in the iteration.
Using this predict-and-recompute approach gives \PipePRMCG and \PipePRCG.

\begin{algorithm}
\caption{Pipelined predict-and-recompute conjugate gradient}\label{alg:pipe_pr_pcg}
\fontsize{10}{10}\selectfont
\begin{algorithmic}[1]
\Procedure{\PipePRCG}{$\vec{A}$, $\vec{M}$, $\vec{b}$, $\vec{x}_0$}
\State \INIT
\For {\( k=1,2,\ldots \)}
    \State \( \vec{x}_k = \vec{x}_{k-1} + \alpha_{k-1} \vec{p}_{k-1} \) \label{pipe_pr_pcg:xk}
    \State \( \vec{r}_k = \vec{r}_{k-1} - \alpha_{k-1} \vec{s}_{k-1} \) \label{pipe_pr_pcg:rk}
    ,~  \( \tilde{\vec{r}}_k = \tilde{\vec{r}}_{k-1} - \alpha_{k-1} \tilde{\vec{s}}_{k-1} \)
    \State \( \vec{w}'_k = \vec{w}_{k-1} - \alpha_{k-1} \vec{u}_{k-1} \)%
    ,~  \( \tilde{\vec{w}}'_k = \tilde{\vec{w}}_{k-1} - \alpha_{k-1} \tilde{\vec{u}}_{k-1} \) \label{pipe_pr_pcg:wkp}
    \State \( \nu'_{k} = \nu_{k-1} - 2 \alpha_{k-1}\sigma_{k-1} + \alpha_{k-1}^2 \gamma_{k-1}  \) \label{pipe_pr_pcg:nukp}
    \State \( \beta_{k} = \nu'_{k} / \nu_{k-1} \) \label{pipe_pr_pcg:betak}
    \State \( \vec{p}_k = \tilde{\vec{r}}_k + \beta_k \vec{p}_{k-1} \) \label{pipe_pr_pcg:pk}
    \State \( \vec{s}_k = \vec{w}'_k + \beta_k \vec{s}_{k-1} \)%
    ,~  \( \tilde{\vec{s}}_k = \tilde{\vec{w}}'_k + \beta_k \tilde{\vec{s}}_{k-1} \) \label{pipe_pr_pcg:sk}
    \State \( \vec{u}_k = \vec{A} \tilde{\vec{s}}_k \)%
    ,~  \( \tilde{\vec{u}}_k = \vec{M}^{-1} \vec{u}_k \) \label{pipe_pr_pcg:uk}
    \State \( \vec{w}_k = \vec{A}\tilde{\vec{r}}_k \)%
    ,~ \( \tilde{\vec{w}}_k = \vec{M}^{-1} \vec{w}_k \) \label{pipe_pr_pcg:wk}
    \State \( \mu_k = \langle \vec{p}_k, \vec{s}_k \rangle \)%
    ,~ \( \sigma_k = \langle \tilde{\vec{r}}_k, \vec{s}_k\rangle \)%
    ,~ \( \gamma_k = \langle \tilde{\vec{s}}_k, \vec{s}_k \rangle \)%
    ,~ \( \nu_k = \langle \tilde{\vec{r}}_k, \vec{r}_k\rangle \) \label{pipe_pr_pcg:inners}
    \State \( \alpha_k = \nu_k / \mu_k \) \label{pipe_pr_pcg:alphak}
\EndFor
\EndProcedure
\end{algorithmic}
\end{algorithm}

\cref{alg:pipe_pr_pcg} shows \PipePRCG, from which \PipePRMCG can be obtained by using the alternate expression \( \nu'_k = -\nu_{k-1} + \alpha_{k-1}^2 \gamma_{k-1}  \) in line~\ref{pipe_pr_pcg:nukp}.
Using this expression means that \( \sigma_k \) need not be computed.
As before, we use a prime to denote predicted quantities.


\subsubsection{Implementation}
\label{sec:implementation}

The presentation of \PipePRCG in \cref{alg:pipe_pr_pcg} is intended to match the derivation from \HSCG and to emphasize the mathematical equivalence of the two algorithms.
However, as with any parallel algorithm, some care must be taken at implementation time as an inefficient implementation may actually increase the runtime per iteration.

We suggest computing the scalars \( \alpha_{k-1}, \nu_k', \) and \( \beta_k \) (lines \ref{pipe_pr_pcg:alphak}, \ref{pipe_pr_pcg:nukp}, \ref{pipe_pr_pcg:betak}) at the beginning of each iteration. 
This will allow all vector updates (lines \ref{pipe_pr_pcg:xk}, \ref{pipe_pr_pcg:rk}, \ref{pipe_pr_pcg:wkp}, \ref{pipe_pr_pcg:pk}, \ref{pipe_pr_pcg:sk}) to occur simultaneously. 
The vector updates require only local on-node communication and are therefore assumed to be very fast.
Finally, the matrix-vector products/preconditioning (lines \ref{pipe_pr_pcg:uk}, \ref{pipe_pr_pcg:wk}) and inner products (line \ref{pipe_pr_pcg:inners}) can all be computed simultaneously.
As a result, the dominant cost per iteration will be either the time for the global reduction associated with the inner products or with the matrix-vector products, thus giving the runtime \( \max( \Cgr, \Ttmv + \Cmv ) \) as listed in \cref{table:comparison}.

The matrix-vector products (and preconditioning) in lines~\ref{pipe_pr_pcg:uk} and \ref{pipe_pr_pcg:wk} can be computed together using efficient kernels.
In particular, this means that within a single node, \PipePRCG still requires only one pass over \( \vec{A} \) (and \( \vec{M}^{-1} \)) in each iteration.
This is an especially important consideration if \( \vec{A} \) is too large to store in fast memory.
Similarly, three of the inner products involve \( \vec{s}_k \), so the number of passes over \( \vec{s}_k \) can be reduced from three to one.
However, this is likely not to have a noticeable effect until the cost of reading \( \vec{s}_k \) from memory becomes large compared to the reduction time.
Finally, there is no need to store \( \vec{w}_k \) and \( \vec{w}_k' \) as separate vectors.

\section{Rounding error analysis}
\label{sec:rounding}

We give a rounding error analysis which provides some insight into how predict-and-recompute schemes may lead to the improved maximal accuracy and convergence behavior observed on test problems.

The maximal attainable accuracy of a CG algorithm in finite precision is typically analyzed in terms of the residual gap \( \Delta_{\vec{r}_k} := (\vec{b} - \vec{A}\vec{x}_k) - \vec{r}_k \) \cite{greenbaum_97a,sleijpen_vandervorst_fokkema_94}, an expression introduced by Greenbaum in \cite[Theorem 2]{greenbaum_89}.
For many variants, such as \HSCG and those introduced in this paper, it is observed experimentally that the norm of the updated residual \( \vec{r}_k \) decreases to much lower than the machine precision.
As a result, the size of the residual gap \( \Delta_{\vec{r}_k} \) can be used to estimate of the size of the smallest true residual which can be attained in finite precision, thereby giving an estimate of the accuracy of the iterate \( \vec{x}_k \).
Similar analyses have been done for a three-term CG variant \cite{gutknecht_strakos_00}, as well as for \CGCG, \GVCG, and other pipelined CG variants \cite{cools_yetkin_agullo_giraud_vanroose_18,carson_rozloznik_strakos_tichy_tuma_18}.
However, for some variants such as \GVCG, the updated residual \( \vec{r}_k \) may not decrease to well below machine precision, so some care must be taken when interpreting such results.

There is also existing theory about the rate of convergence of a CG implementation in finite precision due to Greenbaum \cite{greenbaum_89}, which extends the work of Paige \cite{paige_71,paige_80}.
This analysis applies to a perturbed Lanczos recurrence, which, in effect, shows that the error norms of a CG algorithm run in finite precision for \( k \) steps correspond to the error norms of exact arithmetic CG applied to a larger matrix whose eigenvalues lie in small intervals about the eigenvalues of \( \vec{A} \), provided that the updated residuals \( \vec{r}_k \) satisfy certain conditions.
This provides a method for applying the well-understood theory about the convergence of CG in exact arithmetic to the finite precision setting.

The conditions required for the analysis in \cite{greenbaum_89} are that (i) successive updated residuals are approximately orthogonal; i.e., \( \langle \vec{r}_k, \vec{r}_{k-1} \rangle \approx 0 \), and (ii) that they approximately satisfy the three-term Lanzos recurrence; i.e., 
\begin{align}
    \label{eqn:lanczos}
    \vec{A}\vec{q}_{k}
    &\approx
    \frac{1}{\alpha_{k-1}} \frac{\|\vec{r}_k\|}{\|\vec{r}_{k-1}\|} \vec{q}_{k+1} + \left( \frac{1}{\alpha_{k-1}} + \frac{\beta_{k-1}}{\alpha_{k-2}} \right) \vec{q}_{k} + \frac{1}{\alpha_{k-2}} \frac{\|\vec{r}_{k-1}\|}{\|\vec{r}_{k-2}\|} \vec{q}_{k-1},
\end{align}
where \( \alpha_k \) and \( \beta_k \) are computed in finite precision and \( \vec{q}_{k+1} := (-1)^k \vec{r}_k / \| \vec{r}_k \| \) is obtained from the updated residuals.

If these conditions are satisfied, then for some small \( \eta \) which depends on the machine precision and on the technical definition of ``approximately'' (see \cite{greenbaum_89} for details), the error norms will satisfy the relaxed minimax bound\footnote{%
Greenbaum's theory actually applies to \( \| \vec{A}^{-1/2} \vec{r}_k \| \), which in exact arithmetic is equal to the \( \vec{A} \)-norm of the error \( \| \vec{e}_k \|_{\vec{A}} = \| \vec{A}^{-1} \vec{b} - \vec{x}_k \|_{\vec{A}} = \| \vec{A}^{-1/2} ( \vec{b} - \vec{A} \vec{x}_k ) \| \).
Thus, the relaxed minimax bound \cref{eqn:relaxed_minimax} only applies before the updated and true residuals begin to differ significantly; i.e., when the residual gap \( \Delta_{\vec{r}_k} \) is still small.
}
\begin{align}
    \frac{\| \vec{e}_k \|_\vec{A}}{\| \vec{e}_0 \|_\vec{A}} \leq \min_{p\in\mathcal{P}_k} \left[ \max_{z\in \mathcal{L}(\vec{A})} |p(z)| \right]
    ,&& \mathcal{L}(\vec{A}) = \bigcup_{i=1}^{n} [\lambda_i - \eta, \lambda_i + \eta].
    \label{eqn:relaxed_minimax}
\end{align}
The degree to which these conditions are satisfied in finite precision directly impacts the size of \( \eta \), with better approximations yielding smaller \( \eta \) and therefore stronger bounds.
In fact, if the conditions are exactly satisfied, then the analysis will yield \( \eta = 0 \) and \cref{eqn:relaxed_minimax} becomes the well-known minimax bound for exact arithmetic CG.

However, it remains to be proved that \emph{any} of the variants discussed in this paper satisfy both of the conditions of \cite{greenbaum_89} in a meaningful way.
Experimentally, it seems to be the case that these variants do keep successive residuals approximately orthogonal (this has been essentially proved for \HSCG \cite[Proposition 5.19]{meurant_06}), and that, with the exception of \GVCG, they approximately satisfy the three-term Lanczos recurrence.
In \cite{greenbaum_liu_chen_19}, expressions for the degree to which the updated residuals from \HSCG, \CGCG, and \GVCG satisfy \cref{eqn:lanczos} in finite precision are given in terms of roundoff errors, and it is shown that while \HSCG and \CGCG satisfy the three-term recurrence to within local rounding errors, \GVCG does not.
While this analysis does not prove the degree to which any of the variants satisfy the conditions of \cite{greenbaum_89}, it does provide some indication that the size of \( \eta \) for \HSCG and \CGCG will likely be smaller than for \GVCG on most problems.

Experimentally, on many problems the rate of convergence of all variants is roughly the same prior to the stagnation of the true residuals.
However, on other problems, the rate of convergence of different variants is observed to differ significantly.
Both cases are shown in the experiments in \cref{sec:numerical}; see, e.g., \cref{fig:error_A_norm}. 
In \cite{greenbaum_liu_chen_19} it is suggested that the first case corresponds to problems where the size of \( \eta \) has little effect on the strength of the relaxed minimax bound \cref{eqn:relaxed_minimax} (such as those with a relatively uniformly distributed spectrum), while the second case corresponds to problems where the size of \( \eta \) has a large effect on the strength of this bound (such as those with large gaps in the upper spectrum).
Thus, on the second type of problem where the size of \( \eta \) is important, the fact that \GVCG does not do a good job of satisfying the three-term Lanczos recurrence \cref{eqn:lanczos} means that the relaxed minimax bound will be stronger for \HSCG and \CGCG than for \GVCG.
Again, while this analysis does not prove that \HSCG and \CGCG will have better rates of convergence on such problems, it is evidence supporting the observation that they do have better rates of convergence in practice.

Throughout the analysis we assume the standard model of floating point arithmetic, i.e., 
\begin{align}
    \label{eqn:fp}
    | \operatorname{fp}( \alpha \circ \beta ) - \alpha \circ \beta |
    \leq \epsilon | \alpha \circ \beta |
\end{align}
for real numbers \( \alpha \) and \( \beta \) with standard operations \( \circ \in \{ +, -, \times, \div \} \), where \( \epsilon \) is the unit roundoff of the machine and all computations inside the \( \operatorname{fp}(\cdot) \) are performed in finite precision.
From this, to first order we have the bounds  
\begin{align}
    \label{eqn:axpy}
    \| \operatorname{fp}(\vec{x} + \alpha \vec{y}) - (\vec{x} + \alpha \vec{y}) \| 
    &\leq \epsilon \, ( \|\vec{x}\| + 2 |\alpha | \|\vec{y}\| ), \\
    \label{eqn:Ax}
    \| \operatorname{fp}(\vec{A} \vec{x}) - \vec{A} \vec{x} \| 
    &\leq \epsilon \, c \: \|\vec{A}\| \|\vec{x}\|, \\
    \label{eqn:ip}
    \| \operatorname{fp}(\langle \vec{x}, \vec{y} \rangle ) - \langle \vec{x}, \vec{y} \rangle \| 
    &\leq \epsilon \, n \, \|\vec{x}\| \|\vec{y}\|,
\end{align}
where \( n \) is the length of the vectors \( \vec{x} \) and \( \vec{y} \) and \( c \) is a constant depending on sparsity/structure of \( \vec{A} \) and the method of matrix multiplication used; for instance, it is common to take \( c = m n^{1/2} \) where \( m \) is the maximum number of nonzero entries in a row of \( \vec{A} \).

\subsection{Predict-and-recompute CG}

In finite precision, \PRCG generates the recurrences
\begin{align}
		\begin{aligned}[t]
		\vec{x}_k &= \vec{x}_{k-1} + \alpha_{k-1} \vec{p}_{k-1} + \delta_{\vec{p}_{k}}, \\
		 \nu_k' &= \nu_{k-1} - 2 \alpha_{k-1} \sigma_{k-1} + \alpha_{k-1}^2 \gamma_{k-1} + \delta_{\nu_k'},\\
		\vec{p}_{k} &= \vec{r}_k + \beta_k \vec{p}_{k-1} + \delta_{\vec{p}_{k}},\\
		\nu_k &= \langle \vec{r}_k, \vec{r}_k \rangle + \delta_{\nu_k},\\
		\gamma_k &= \langle \vec{s}_k, \vec{s}_k \rangle + \delta_{\gamma_k},
		\end{aligned}
\qquad
		\begin{aligned}[t]
		\vec{r}_k &= \vec{r}_{k-1} - \alpha_{k-1} \vec{s}_{k-1} + \delta_{\vec{r}_k},\\
	    \beta_k &= \nu_k' / \nu_{k-1} + \delta_{\beta_k},\\
		\vec{s}_{k} &= \vec{A} \vec{p}_{k} + \delta_{\vec{s}_{k}},\\
		\sigma_k &= \langle \vec{r}_k, \vec{s}_k \rangle + \delta_{\sigma_k},\\
	\end{aligned}
		\label{eqn:PRerrors}
\end{align}
where the \( \delta_{\text{quantity}} \) terms represent roundoff errors incurred at each step and are bounded by the appropriate application of \cref{eqn:fp,eqn:axpy,eqn:Ax,eqn:ip}.

It was observed by Meurant~\cite{meurant_87} that the instability caused by recursive computation of \( \nu_k \) via the formula \cref{eqn:nuk_simplified} is due to the value of \( \nu_k \) eventually becoming negative. Our methods, like those in \cite{rosendale_83,saad_85,meurant_87}, also break down if \( \nu_k \) becomes negative. Meurant conjectures in~\cite{meurant_87} that using the predict-and-recompute approach appears to solve this problem, but he does not offer a rigorous theoretical explanation of this behavior. 

We begin this section by filling in this gap. We first show that the predict-and-recompute scheme suggested in \cite{meurant_87} and used in \MCG, \PRCG, \PipePRMCG, and \PipePRCG keeps the estimate for \( \nu_k' \) to within local rounding errors of the true value \( \langle \vec{r}_k, \vec{r}_k \rangle \).
To this end, we define the \( \nu' \)-gap as \( \Delta_{\nu_k'} := \langle \vec{r}_k, \vec{r}_k \rangle - \nu_k' \).
Then
\begin{align}
    \Delta_{\nu_k'}
    \tag*{} &= 
    \langle \vec{r}_{k-1} - \alpha_{k-1} \vec{s}_{k-1} + \delta_{\vec{r}_k}, \vec{r}_{k-1} - \alpha_{k-1} \vec{s}_{k-1} + \delta_{\vec{r}_k} \rangle - \nu_k' 
    \\ \tag*{} &= 
    \langle \vec{r}_{k-1}, \vec{r}_{k-1} \rangle - 2 \alpha_{k-1} \langle \vec{r}_{k-1}, \vec{s}_{k-1} \rangle + \alpha_{k-1}^2 \langle \vec{s}_{k-1}, \vec{s}_{k-1} \rangle  
    \\ \tag*{}&\hspace{10em} 
    + 2 \langle \delta_{\vec{r}_k}, \vec{r}_{k-1} - \alpha_{k-1} \vec{s}_{k-1} \rangle + \langle \delta_{\vec{r}_k}, \delta_{\vec{r}_k} \rangle - \nu_k' 
    \\ \tag*{} &= 
		\nu_{k-1} - \delta_{\nu_{k-1}} - 2\alpha_{k-1}(\sigma_{k-1} - \delta_{\sigma_{k-1}}) + \alpha_{k-1}^2 (\gamma_{k-1} - \delta_{\gamma_{k-1}}) 
			\\ \tag*{}&\hspace{10em}
			+ 2 \langle \delta_{\vec{r}_k}, \vec{r}_{k} \rangle - \langle \delta_{\vec{r}_k}, \delta_{\vec{r}_k} \rangle - \nu_k' 
		\\ \tag*{} &= (\nu_{k-1} - 2 \alpha_{k-1} \sigma_{k-1} + \alpha_{k-1}^2 \gamma_{k-1}) - (\delta_{\nu_{k-1}} - 2 \alpha_{k-1} \delta_{\sigma_{k-1}} + 		\alpha_{k-1}^2 \delta_{\gamma_{k-1}}) 
			\\ \tag*{}&\hspace{10em} 
			+ 2 \langle \delta_{\vec{r}_k}, \vec{r}_{k} \rangle - \langle \delta_{\vec{r}_k}, \delta_{\vec{r}_k}		\rangle - \nu_k' 
		\\ \tag*{} &= \nu_{k}' - \delta_{\nu_{k}'} - (\delta_{\nu_{k-1}} - 2 \alpha_{k-1} \delta_{\sigma_{k-1}} + \alpha_{k-1}^2 \delta_{\gamma_{k-1}})
		+  2 \langle \delta_{\vec{r}_k}, \vec{r}_{k} \rangle - \langle \delta_{\vec{r}_k}, \delta_{\vec{r}_k}	\rangle - \nu_k'
		\\ &= 2\langle \delta_{\vec{r}_k}, \vec{r}_{k} \rangle - \langle \delta_{\vec{r}_k}, \delta_{\vec{r}_k}	\rangle - \delta_{\nu_{k}'} 
		- (\delta_{\nu_{k-1}} - 2 \alpha_{k-1} \delta_{\sigma_{k-1}} + \alpha_{k-1}^2 \delta_{\gamma_{k-1}}) \label{eqn:Delta_nukp}.
\end{align}
Applying \cref{eqn:axpy,eqn:Ax,eqn:ip} we have the bound
\begin{align}
    | \Delta_{\nu_k'} |
    \leq 2 \| \delta_{\vec{r}_k}\| \| \vec{r}_k \| 
    + \| \delta_{\vec{r}_k} \|^2 
    + | \delta_{\nu'_k} | 
    + | \delta_{\nu_{k-1}} | 
    + 2 | \alpha_{k-1} | | \delta_{\sigma_{k-1}} | + |\alpha_{k-1}|^2 | \delta_{\gamma_{k-1}} |.
\
\label{eqn:nuk_bound}
\end{align}

We will now bound the terms on the right-hand side of~\cref{eqn:nuk_bound} using standard analysis techniques. 
We can rearrange the expression for \( \vec{r}_{k} \) in \cref{eqn:PRerrors} to obtain 
\begin{equation*}
    \vec{s}_{k-1} 
    = \frac{1}{|\alpha_{k-1}|}(\vec{r}_{k-1}-\vec{r}_k + \delta_{\vec{r}_k}),
\end{equation*}
which gives the bound 
\begin{equation}
    \| \vec{s}_{k-1} \| 
    \leq \frac{1}{|\alpha_{k-1}|} \left( \| \vec{r}_{k-1}\| + \| \vec{r}_k \| + \| \delta_{\vec{r}_k}\| \right).
\label{eqn:sbound}
\end{equation}

\begin{figure*}[t]\centering
    \begin{subfigure}{.48\textwidth}
        \includegraphics[width=\textwidth]{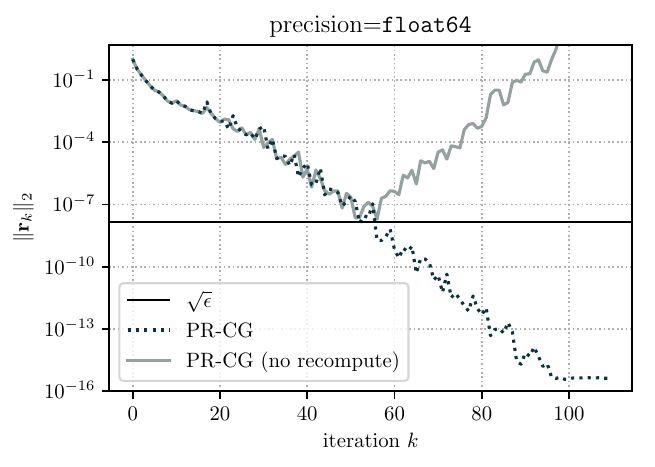}
    \end{subfigure}
    \begin{subfigure}{.48\textwidth}
        \includegraphics[width=\textwidth]{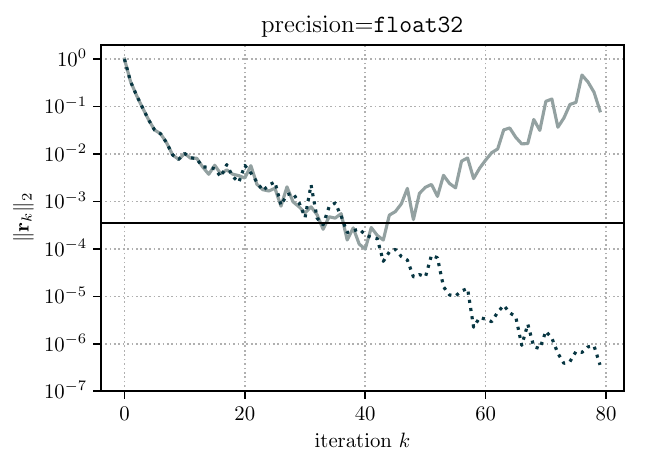}
    \end{subfigure}
    \caption{Demonstration, on \texttt{model\_48\_8\_3} with no preconditioner, of breakdown of \PRCG when \( \nu_k' \) is not recomputed.
    Breakdown occurs when \( \nu_k' \) becomes negative which can be avoided using the recompute scheme.
    }
    \label{fig:no_recompute}
\end{figure*}

We now seek to bound the norms of the quantities \(\delta_{\vec{r}_k}\), \(\delta_{\nu_{k-1}}\), \(\delta_{\sigma_{k-1}}\), \(\delta_{\gamma_{k-1}}\), and \(\delta_{\nu'_k}\). Using \cref{eqn:PRerrors}, \cref{eqn:sbound}, and the bounds \cref{eqn:axpy}-\cref{eqn:ip}, and dropping terms of order \(O(\epsilon^2)\), we have 
\begin{align}
    \| \delta_{\vec{r}_k} \| 
    &\leq \epsilon \left( \| \vec{r}_k \| + 2 |\alpha_{k-1} \| \vec{s}_{k-1} \| \right)\nonumber \\
    &\leq \epsilon \left( \| \vec{r}_k \| + 2 |\alpha_{k-1} | \left( \frac{1}{|\alpha_{k-1}|} \left( \| \vec{r}_{k-1}\| + \| \vec{r}_k \| + \| \delta_{\vec{r}_k}\| \right) \right) \right) \nonumber \\
    &\leq \epsilon \left( 3 \| \vec{r}_{k}\| + 2 \| \vec{r}_{k-1}\| \right) + 2\,\epsilon \, \| \delta_{\vec{r}_{k}}\| \nonumber \\
    &\leq 3\, \epsilon \left( \| \vec{r}_k \| + \| \vec{r}_{k-1} \| \right), 
    \label{eqn:drbound} \\
    | \delta_{\nu_{k-1}} | 
    &\leq \epsilon\, n\, \| \vec{r}_{k-1}\|^2, \label{eqn:dnubound}\\
    | \delta_{\sigma_{k-1}} | 
    &\leq \epsilon\, n\, \| \vec{r}_{k-1}\| \|\vec{s}_{k-1}\| \nonumber \\
    &\leq \epsilon\, n\, \| \vec{r}_{k-1}\| \left(\frac{1}{|\alpha_{k-1}|} \left( \|\vec{r}_{k-1}\| + \|\vec{r}_k\| + \|\delta_{\vec{r}_k}\| \right) \right) \nonumber \\
    &\leq \epsilon\, n\, \frac{1}{|\alpha_{k-1}|} \left( \|\vec{r}_{k-1}\| ^2 + \| \vec{r}_{k-1} \| \| \vec{r}_k \| \right) , \label{eqn:ddbound}\\
    | \delta_{\gamma_{k-1}} | 
    &\leq \epsilon\, n\, \| \vec{s}_{k-1}\|^2 \nonumber \\
    &\leq \epsilon\, n\, \frac{1}{|\alpha_{k-1}|^2} \left( \| \vec{r}_{k-1}\| + \| \vec{r}_{k}\| \right)^2 ,  \label{eqn:dgbound}
\intertext{
where we have used the bound,
\begin{align*}
    \| \vec{s}_{k-1} \| \leq (1 + 3 \epsilon) \frac{1}{|\alpha_{k-1}|}(\|\vec{r}_{k-1}\| + \|\vec{r}_{k} \|).
\end{align*}
Finally, using \cref{eqn:PRerrors}, \cref{eqn:sum_prod}, and \cref{eqn:axpy,eqn:Ax,eqn:ip}, we find
}
    | \delta_{\nu_k'} |
    &\leq 3 \, \epsilon \left( | \nu_{k-1} | + 2 | \alpha_{k-1} | | \sigma_{k-1} | + |\alpha_{k-1}|^2 | \gamma_{k-1} | \right) \nonumber
    \\&\leq 3\, \epsilon \left( \| \vec{r}_{k-1} \|^2 
        + 2 | \alpha_{k-1} | \| \vec{r}_{k-1} \| \| \vec{s}_{k-1} \|
            + |\alpha_{k-1}|^2 \| \vec{s}_{k-1} \|^2 \right) \nonumber
    \\&\leq 3\, \epsilon \left( \| \vec{r}_{k-1} \|^2 
        + 2 \| \vec{r}_{k-1} \| (\| \vec{r}_{k-1} \| 
        + \| \vec{r}_k \|) 
        + (\| \vec{r}_{k-1} \| + \| \vec{r}_k \|)^2  \right) \nonumber
    \\&= 3\, \epsilon \left( 4 \| \vec{r}_{k-1} \|^2 + 4\| \vec{r}_{k-1} \|\| \vec{r}_{k} \| + \| \vec{r}_{k} \|^2  \right)
    \label{eqn:dnpbound}.
\end{align}

Substituting \cref{eqn:drbound,eqn:dnubound,eqn:ddbound,eqn:dgbound,eqn:dnpbound} into~\cref{eqn:nuk_bound} gives a bound on the \( \nu_k' \) gap:
\begin{align}
    | \Delta_{\nu_k'} | 
    &\leq \epsilon \left( (12+3n)\| \vec{r}_{k-1} \|^2 
    + (18+4n) \| \vec{r}_{k-1} \|\| \vec{r}_{k} \| 
    + (9+n) \| \vec{r}_{k} \|^2 \right) 
    \\&\leq \epsilon \, (30+7n) \left( \| \vec{r}_{k-1} \|^2 + \| \vec{r}_{k} \|^2 \right) .
\label{eqn:nuk_bound_final}
\end{align}
Therefore, given that \(\nu'_k = \| \vec{r}_k \|^2 - \Delta_{\nu'_k}\), it is clear from the above bound that 
\begin{equation*}
    \nu'_k 
    \geq \| \vec{r}_k \|^2 - \epsilon \, (30+7n) \left( \| \vec{r}_{k-1} \|^2 + \| \vec{r}_{k} \|^2 \right) 
\end{equation*}
Assuming \( \epsilon \) is small enough that the higher order terms do not have an significant impact on the derived bounds, this quantity will remain positive provided \( \| \vec{r}_{k-1} \|^2 / \| \vec{r}_k \|^2 \leq 1 / (\epsilon \, (30 + 7 n)) - 1 \).
While we do not explicitly bound \( \| \vec{r}_{k-1} \| / \| \vec{r}_k \| \), we note that in practice it is unlikely to be of size \( \sim 1/\sqrt{\epsilon n} \), thus explaining why the recompute scheme keeps \( \nu_k' \) positive.

We now consider the case where we do not use the recomputed value \( \nu_{k-1} \), i.e., if we were to instead compute 
\begin{align} 
    \nu_k' = \nu_{k-1}' - 2 \alpha_{k-1} \sigma_{k-1} + \alpha_{k-1}^2 \gamma_{k-1} + \underline{\delta}_{\nu_{k}'},
    \label{eqn:nuk_norecomp}
\end{align}
where \(\underline{\delta}_{\nu_{k}'}\) represents the local rounding errors in computing this recurrence.  
In this case we would have, letting \( \underline\Delta_{\nu_k'} \) represent the \( \nu' \)-gap without recomputation, 
\begin{align*}
    \underline\Delta_{\nu_k'} 
    &= \langle \vec{r}_{k-1}, \vec{r}_{k-1} \rangle - 2\alpha_{k-1} (\sigma_{k-1} - \delta_{\sigma_{k-1}}) + \alpha_{k-1}^2 (\gamma_{k-1} - \delta_{\gamma_{k-1}}) \\
	&\phantom{=} + 2 \langle \delta_{\vec{r}_k}, \vec{r}_k \rangle - \langle \delta_{\vec{r}_k}, \delta_{\vec{r}_k} \rangle - \nu_{k-1}' + 2\alpha_{k-1}\sigma_{k-1} - \alpha_{k-1}^2 \gamma_{k-1} - \underline{\delta}_{\nu_k'} \\
    &= (\langle \vec{r}_{k-1}, \vec{r}_{k-1} \rangle - \nu_{k-1}') + 2\alpha_{k-1}\delta_{\sigma_{k-1}} - \alpha_{k-1}^2 \delta_{\gamma_{k-1}} 
    +2 \langle \delta_{\vec{r}_k}, \vec{r}_k \rangle - \langle \delta_{\vec{r}_k}, \delta_{\vec{r}_k} \rangle  - \underline{\delta}_{\nu_k'}\\
    &= \underline\Delta_{\nu_{k-1}'}+ 2\alpha_{k-1}\delta_{\sigma_{k-1}} - \alpha_{k-1}^2 \delta_{\gamma_{k-1}} 
    +2 \langle \delta_{\vec{r}_k}, \vec{r}_k \rangle - \langle \delta_{\vec{r}_k}, \delta_{\vec{r}_k} \rangle  - \underline{\delta}_{\nu_k'}.
%
\end{align*}
Taking the absolute value of both sizes and and ignoring terms of \( \mathcal{O}(\epsilon) \), we thus have the bound 
\begin{equation}
    | \underline\Delta_{\nu_k'} | 
    \leq | \underline\Delta_{\nu'_{k-1}} | 
        + 2 |\alpha_{k-1}| | \delta_{\sigma_{k-1}} | 
        + \alpha_{k-1}^2 | \delta_{\gamma_{k-1}} | 
        + 2 \| \delta_{\vec{r}_k}\| \|\vec{r}_k\| 
        + | \underline{\delta}_{\nu_k'} |.
\label{eqn:Dunuk_bound}
\end{equation}

The only term in \cref{eqn:Dunuk_bound} that remains to bound is \( | \underline{\delta}_{\nu_k'} | \).
From \cref{eqn:nuk_norecomp} we have
\[
    | \underline{\delta}_{\nu_k'} | 
    \leq 3\,\epsilon \left( | \nu'_{k-1} | + 2 |\alpha_{k-1}| | \sigma_{k-1} | + \alpha_{k-1}^2 \|\gamma_{k-1}\| \right).
\]
By definition, \(\nu'_{k-1}=\| \vec{r}_{k-1}\|^2 -\underline{\Delta}_{\nu'_{k-1}} \), so
\begin{align*}
    | \underline{\delta}_{\nu_k'} | 
    &\leq 3\, \epsilon \left( \| \vec{r}_{k-1}\|^2 + | \underline{\Delta}_{\nu'_{k-1}} | + 2 |\alpha_{k-1}| | \sigma_{k-1} | + \alpha_{k-1}^2 | \gamma_{k-1} | \right).
\end{align*}
%
Thus, using an approach similar to the computation in \cref{eqn:dnpbound} we have
\begin{equation}
    | \underline{\delta}_{\nu_k'} | 
    \leq 3\, \epsilon \left( 4 \| \vec{r}_{k-1} \|^2 + 4\| \vec{r}_{k-1} \|\| \vec{r}_{k} \| + \| \vec{r}_{k} \|^2 + | \underline{\Delta}_{\nu'_{k-1}} | \right).
\label{eqn:dunuk_bound}
\end{equation}
Plugging \cref{eqn:dunuk_bound}, \cref{eqn:ddbound}, \cref{eqn:dgbound}, and \cref{eqn:drbound} into \cref{eqn:Dunuk_bound}, we obtain 
\begin{align*}
    | \underline\Delta_{\nu_k'} | 
    \leq (1 + 3 \, \epsilon) | \underline\Delta_{\nu_{k-1}'} | + \epsilon\, \mathcal{O}(n) \left( \| \vec{r}_{k-1} \|^2 + \| \vec{r}_k \|^2 \right).
\end{align*}
Using that \( | \underline{\Delta}_{\nu'_{0}} | \leq \epsilon \, n \, \|\vec{r}_0\|^2\), this can be written 
\begin{equation}
    | \underline\Delta_{\nu_k'} |
    \leq \epsilon\,\mathcal{O}(n) \sum_{i=0}^{k} \| \vec{r}_{i} \|^2,
\label{eqn:Dunuk_bound_final}
\end{equation}
Thus for the $\nu'$-gap without recomputation, we have the bound
\begin{align}
    \nu'_k \geq \|\vec{r}_k\|^2 - \epsilon\,\mathcal{O}(n) \sum_{i=0}^{k} \| \vec{r}_{i} \|^2.
    \label{eqn:nuk_no_recompute}
\end{align}

Since now the term subtracted involves the sum of the squares of \emph{all} previous residuals, which can be large especially at the  beginning of the iterations, it is entirely possible that \(\nu'_k\) can become negative at some point during the iterations.
While \cref{eqn:nuk_no_recompute} is only an upper bound, it is more or less clear that this is the cause of instability in the case without recomputation, which explains the observations of Meurant~\cite{meurant_87}.
In fact, if we assume that early residuals have norm \( \approx 1 \), the form of the expression \cref{eqn:nuk_no_recompute} suggests that \( \nu_k' \) can no longer be guaranteed to be positive once \( \| \vec{r}_k \| \approx \sqrt{ \epsilon } \).
This aligns with what we have observed in practice; see for instance \cref{fig:no_recompute}.

\subsubsection{Maximum attainable accuracy}

We now compute an expression for the residual gap in \cref{alg:pr_pcg}.
Substituting in our recurrences for \( \vec{x}_k \) and \( \vec{r}_k \) computed in finite precision, we find
\begin{align}
    \Delta_{\vec{r}_k} 
    \tag*{} &=
    \vec{b} - \vec{A} (\vec{x}_{k-1} + \alpha_{k-1} \vec{p}_{k-1} + \delta_{\vec{x}_k}) - (\vec{r}_{k-1} - \alpha_{k-1} \vec{s}_{k-1} + \delta_{\vec{r}_k})
    \\ \tag*{} &=
    ( \vec{b} - \vec{A}\vec{x}_{k-1} - \vec{r}_{k-1}) + \alpha_{k-1}(\vec{s}_{k-1} - A \vec{p}_{k-1}) - \vec{A}\delta_{\vec{x}_k} - \delta_{\vec{r}_k}
    \\ \tag*{} &=\Delta_{\vec{r}_{k-1}} + \alpha_{k-1} \delta_{\vec{s}_{k-1}} - \vec{A} \delta_{\vec{x}_k} - \delta_{\vec{r}_k}.
\end{align}
Therefore the residual gap can be written 
\begin{align*}
    \Delta_{\vec{r}_k}
    = \Delta_{\vec{r}_0} + \sum_{i=1}^{k} \left( \alpha_{i-1} \delta_{\vec{s}_{i-1}} - \vec{A}\delta_{\vec{x}_i} - \delta_{\vec{r}_i} \right).
\end{align*}
This expression is a simple accumulation of local rounding errors.
We note that the residual gaps for \HSCG and \CGCG are both given by
\begin{align*}
    \Delta_{\vec{r}_k}
    = \Delta_{\vec{r}_0} - \sum_{i=1}^{k} \left( \vec{A}\delta_{\vec{x}_i} + \delta_{\vec{r}_i} \right),
\end{align*}
where the \( \delta_{\text{quantity}} \) terms correspond to the round off errors made by those algorithms respectively \cite{cools_yetkin_agullo_giraud_vanroose_18}.

\subsubsection{Towards understanding convergence}

We now derive an expression for the extent to which the recurrence \cref{eqn:lanczos} is satisfied by \PRCG in finite precision. 
We express \( \vec{s}_k = \vec{A} \vec{p}_k + \delta_{\vec{s}_k} \) in terms of \( \vec{r}_k \) and \( \vec{s}_{k-1} \), writing 
\begin{align}
    \vec{s}_{k} 
    \tag*{}&= 
    \vec{A} ( \vec{r}_{k} + \beta_k \vec{p}_{k-1} + \delta_{\vec{p}_{k}}) + \delta_{\vec{s}_{k}}
    \\ \tag*{}&= 
    \vec{A} \vec{r}_k + \beta_k \vec{A} \vec{p}_{k-1} + \vec{A} \delta_{\vec{p}_{k}} + \delta_{\vec{s}_{k}}
    \\ &= 
    \vec{A} \vec{r}_k + \beta_k ( \vec{s}_{k-1} - \delta_{\vec{s}_{k-1}} ) + \vec{A} \delta_{\vec{p}_{k}} + \delta_{\vec{s}_{k}}.
\label{eqn:prcg:sk}
\end{align}
By rearranging the equation \( \vec{r}_{k} = \vec{r}_{k-1} - \alpha_{k-1} \vec{s}_{k-1} + \delta_{\vec{r}_k} \), we can write \( \vec{s}_{k-1} = (1/\alpha_{k-1}) (\vec{r}_{k-1}-\vec{r}_{k} + \delta_{\vec{r}_{k}}) \).
From this we compute
\begin{align*}
    \vec{s}_k - \beta_k \vec{s}_{k-1}
    &= \frac{1}{\alpha_k} \left( \vec{r}_k - \vec{r}_{k+1} + \delta_{\vec{r}_{k+1}} \right) - \frac{\beta_k}{\alpha_{k-1}} \left( \vec{r}_{k-1} - \vec{r}_k + \delta_{\vec{r}_k} \right)
    \\&= - \frac{1}{\alpha_k} \vec{r}_{k+1} + \left( \frac{1}{\alpha_k} + \frac{\beta_k}{\alpha_{k-1}} \right) \vec{r}_k - \frac{\beta_k}{\alpha_{k-1}} \vec{r}_{k-1} + \frac{1}{\alpha_k} \delta_{\vec{r}_{k+1}} - \frac{\beta_k}{\alpha_{k-1}} \delta_{\vec{r}_k}. 
\end{align*}

Thus, shifting the indexing down by one, \cref{eqn:prcg:sk} becomes
\begin{align}
    \vec{A}\vec{r}_{k-1}
    &=
    -\frac{1}{\alpha_{k-1}} \vec{r}_k + \left( \frac{1}{\alpha_{k-1}} + \frac{\beta_{k-1}}{\alpha_{k-2}} \right) \vec{r}_{k-1} - \frac{\beta_{k-1}}{\alpha_{k-2}} \vec{r}_{k-2} + \vec{f}_k,
    \label{eqn:prcg:three_term_r}
\end{align}
where
\begin{align*}
   \vec{f}_k
    &= 
    - \frac{\beta_{k-1}}{\alpha_{k-2}} \delta_{\vec{r}_{k-1}} + \beta_{k-1} \delta_{\vec{s}_{k-2}} - \vec{A} \delta_{\vec{p}_{k-1}} - \delta_{\vec{s}_{k-1}} + \frac{1}{\alpha_{k-1}}\delta_{\vec{r}_k}. 
\end{align*}
Defining \( \vec{q}_{k+1} := (-1)^{k} \vec{r}_k / \| \vec{r}_k \| \) we obtain the three-term recurrence
\begin{align}
    \vec{A}\vec{q}_{k}
    &=
    \frac{1}{\alpha_{k-1}} \frac{\|\vec{r}_k\|}{\|\vec{r}_{k-1}\|} \vec{q}_{k+1} + \left( \frac{1}{\alpha_{k-1}} + \frac{\beta_{k-1}}{\alpha_{k-2}} \right) \vec{q}_{k} + \frac{\beta_{k-1}}{\alpha_{k-2}} \frac{\|\vec{r}_{k-2}\|}{\|\vec{r}_{k-1}\|} \vec{q}_{k-1} - \frac{(-1)^{k-1}}{\|\vec{r}_{k-1}\|} \vec{f}_k.
    \label{eqn:prcg:three_term_q}
\end{align}

In order for Greenbaum's analysis to apply, we must write this as a symmetric three-term recurrence with some perturbation term.
Using our definition of the \( \nu' \)-gap, \( \Delta_{\nu_k'} := \langle  \vec{r}_k , \vec{r}_k \rangle - \nu_k' \), we have
\begin{align*}
    \beta_k
    &=
    \frac{\nu_k'}{\nu_{k-1}} + \delta_{\beta_k}
    = \frac{\|\vec{r}_k\|^2 + \Delta_{\nu_k'}}{\|\vec{r}_{k-1}\|^2 + \delta_{\nu_{k-1}'}} + \delta_{\beta_k}
    = \frac{\|\vec{r}_k\|^2}{\|\vec{r}_{k-1}\|^2} + \Delta_{\beta_k},
\end{align*}
where \( \Delta_{\beta_k} := \beta_k - \| \vec{r}_k \|^2 / \| \vec{r}_{k-1} \|^2 \) can be written explicitly as
\begin{align*}
    \Delta_{\beta_k} 
    = \frac{\Delta_{\nu_k'} + \| \vec{r}_{k-1} \|^2 \delta_{\beta_k} + \delta_{\nu_{k-1}'} \delta_{\beta_k}}{\| \vec{r}_{k-1} \|^2 + \delta_{\nu_{k-1}'}} + \frac{- \| \vec{r}_k \|^2 \delta_{\nu_{k-1}'} }{ \| \vec{r}_{k-1} \|^2(\| \vec{r}_{k-1} \|^2 + \delta_{\nu_{k-1}'})}.
\end{align*}

Then, by plugging this expression for \( \beta_{k-1} \) into \cref{eqn:prcg:three_term_q}, we obtain the approximately symmetric three-term recurrence for \( \vec{q}_k \),
\begin{align*}
    \vec{A}\vec{q}_{k}
    &=
    \frac{1}{\alpha_{k-1}} \frac{\|\vec{r}_k\|}{\|\vec{r}_{k-1}\|} \vec{q}_{k+1} + \left( \frac{1}{\alpha_{k-1}} + \frac{\beta_{k-1}}{\alpha_{k-2}} \right) \vec{q}_{k} + \frac{1}{\alpha_{k-2}} \frac{\|\vec{r}_{k-1}\|}{\|\vec{r}_{k-2}\|} \vec{q}_{k-1} + \vec{F}_{k},
\end{align*}
so that the error in the three-term Lanczos recurrence is given by
\begin{align*}
    \vec{F}_{k}
    &= \frac{\Delta_{\beta_{k-1}}}{\alpha_{k-2}} \frac{\|\vec{r}_{k-2}\|^2}{\|\vec{r}_{k-1}\|^2} \vec{q}_{k-1}  - \frac{(-1)^k}{\|\vec{r}_{k-1}\|} \vec{f}_k.
\end{align*}

Hence, the amount by which the updated residuals in \PRCG fail to satisfy \cref{eqn:lanczos} depends only on local rounding errors, which was also observed for \HSCG and \CGCG in \cite{greenbaum_liu_chen_19}.
While \( \| \vec{r}_{k-1} \| \) may become quite small, we note that \( \vec{f}_k \) consists of roundoff errors made in iterations \( k-1 \) and \( k-2 \) and corresponds to vectors which are of comparable size to \( \vec{r}_{k-1} \); i.e., it is reasonable to assume that the ratio of the norms of these vectors to \( \vec{r}_{k-1} \) will be much less than \( 1/\epsilon \).
Based on this heuristic, \( \vec{f}_k / \| \vec{r}_{k-1} \| \) can be expected to stay small, which provides some indication for why the rate of convergence of \PRCG tends to be similar to \HSCG and \CGCG.

\begin{figure*}[t]\centering
    \begin{subfigure}{.48\textwidth}
        \includegraphics[width=\textwidth]{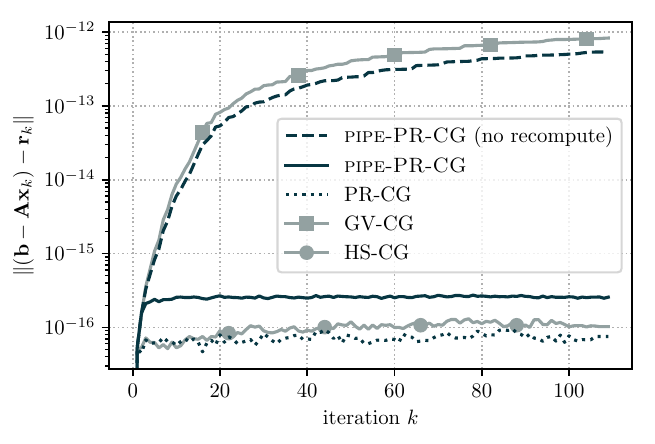}
    \end{subfigure}
    \begin{subfigure}{.48\textwidth}
        \includegraphics[width=\textwidth]{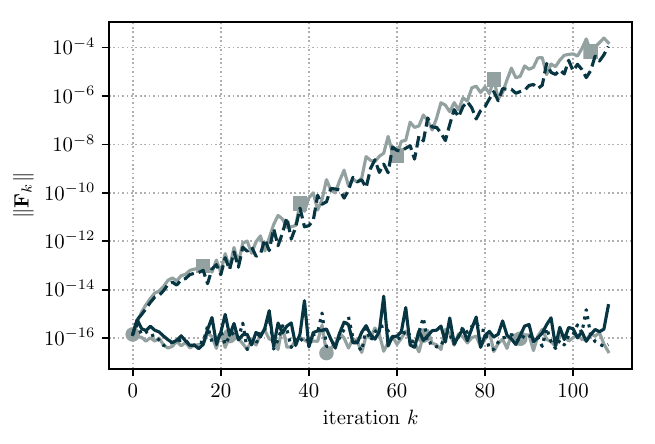}
    \end{subfigure}
    \caption{Residual gap (lefT) and 2-norm of error  in three-term Lanczos recurrence (right) on \texttt{model\_48\_8\_3} without a preconditioner. In the legend, ``no recompute'' indicates that \( \vec{w}_k \) is not recomputed.}
    \label{fig:residual_gap-three_term_error}
\end{figure*}

\subsection{Pipelined predict-and-recompute CG}

We now provide for \PipePRCG a similar rounding error analysis for the residual gap and Lanczos recurrence as we did for \PRCG.
Of course, the \( \delta_{\text{quantity}} \) error terms in this section are different from those in the previous section, even though we use many of the same symbols for notational convenience.
In particular, our expression for \( \vec{s}_k \) in finite precision is now \(\vec{s}_{k} = \vec{w}'_k + \beta_k \vec{s}_{k-1} + \delta_{\vec{s}_{k}} \) and we add the recurrences
\begin{align*}
    \vec{w}'_k &= \vec{w}_{k-1} - \alpha_{k-1} \vec{u}_{k-1} + \delta_{\vec{w}'_k}
    ,&&&
    \vec{w}_k &= \vec{A}\vec{r}_{k} + \delta_{\vec{w}_k}
    ,\\
    \vec{p}_{k} &= \vec{r}_k + \beta_k \vec{p}_{k-1} + \delta_{\vec{p}_{k}}
    ,&&&
    \vec{u}_k &= \vec{A}\vec{s}_{k} + \delta_{\vec{u}_k}. 
\end{align*}
We note that the same expression holds for \( \Delta_{\nu_k'} \) as with \PRCG; see \cref{eqn:Delta_nukp}.

\subsubsection{Maximum attainable accuracy}

Similarly as for \PRCG, we can write the residual gap in \PipePRCG as 
\begin{align*}
    \Delta_{\vec{r}_k} &= \vec{b} - \vec{A}\vec{x}_k - \vec{r}_k \\
		&= \Delta_{\vec{r}_{k-1}} + \alpha_{k-1}(\vec{s}_{k-1} - \vec{A} \vec{p}_{k-1}) - \vec{A} \delta_{\vec{x}_k} - \delta_{\vec{r}_k} \\
		&=\Delta_{\vec{r}_{k-1}} - \alpha_{k-1}\Delta_{\vec{s}_{k-1}} - \vec{A} \delta_{\vec{x}_k} - \delta_{\vec{r}_k} \\
        &=\Delta_{\vec{r}_0} - \sum_{i=1}^{k} \left[ \alpha_{i-1} \Delta_{\vec{s}_{i-1}} + \vec{A}\delta_{\vec{x}_i} + \delta_{\vec{r}_i} \right],
\end{align*}
where \( \Delta_{\vec{s}_i} := \vec{A} \vec{p}_{i} - \vec{s}_i \) represents the gap between \( \vec{A} \vec{p}_{i} \) and the recursively-computed quantity \( \vec{s}_i \). 
We can write this \(\vec{s}\)-gap as 
\begin{align*}
    \Delta_{\vec{s}_{i}} 
    &= 
    \vec{A} ( \vec{r}_i + \beta_{i} \vec{p}_{i-1} + \delta_{\vec{p}_{i}} ) - (\vec{w}_i' + \beta_i \vec{s}_{i-1} + \delta_{\vec{s}_{i}} )
    \\&= 
    (\vec{A} \vec{r}_i- \vec{w}_i') + \beta_i (A \vec{p}_{i-1} - \vec{s}_{i-1}) + \vec{A}\delta_{\vec{p}_{i}} - \delta_{\vec{s}_{i}}
    \\&= 
    \Delta_{\vec{w}_i'} + \beta_i \Delta_{\vec{s}_{i-1}} + \vec{A}\delta_{\vec{p}_{i}} - \delta_{\vec{s}_{i}}\\
		&= \Delta_{\vec{s}_0} \prod_{j=1}^{i} \beta_j + \sum_{j=1}^{i} \Bigg[ (\Delta_{\vec{w}_j'} + \vec{A}\delta_{\vec{p}_j} - \delta_{\vec{s}_j} ) \prod_{\ell=j+1}^{i} \beta_\ell  \Bigg],
\end{align*}
where \( \Delta_{\vec{w}_j'} \) denotes the gap \( \Delta_{\vec{w}_j'} := \vec{A}\vec{r}_j - \vec{w}_j' \). Now we look at how to bound this \( \vec{w'} \)-gap. 
In \PipePRCG, we have 
\begin{align}
    \Delta_{\vec{w}_j'}
    &= 
    \vec{A} (\vec{r}_{j-1} - \alpha_{j-1}\vec{s}_{j-1} + \delta_{\vec{r}_j}) - (\vec{w}_{j-1} - \alpha_{j-1} \vec{u}_{j-1} + \delta_{\vec{w}_j'} ) \nonumber
    \\&= 
    \vec{A} \vec{r}_{j-1} - \vec{w}_{j-1} - \alpha_{j-1} ( \vec{A} \vec{s}_{j-1} - \vec{u}_{j-1}) - \delta_{\vec{w}_k'} + \vec{A} \delta_{\vec{r}_j} \nonumber
    \\&=
    - \delta_{\vec{w}_{j-1}} + \alpha_{j-1} \delta_{\vec{u}_{j-1}} - \vec{A}\delta_{\vec{r}_j} - \delta_{\vec{w}_j'}, \label{eq:Dwr}
\end{align}
which is the sum of local rounding errors. 

Notice that without recomputation, i.e., if we were to instead compute
\begin{equation*}
\vec{w}_j' - \vec{w}_{j-1}' - \alpha_{j-1} \vec{u}_{j-1} + \underline{\delta}_{\vec{w}_j'},
\end{equation*}
where \(\underline{\delta}_{\vec{w}_j'}\) denotes the local rounding errors in computing this recurrence, then 
we would have, letting \( \underline\Delta_{\vec{w}_j'} \) denote the \( \vec{w}' \)-gap in the case of no recomputation, 
\begin{align}
    \underline\Delta_{\vec{w}_j'}
    &=
    \vec{A} \vec{r}_{j-1} - \vec{w}_{j-1}' - \alpha_{j-1} (\vec{A}\vec{s}_{j-1} - \vec{u}_{j-1}) + \vec{A}\delta_{\vec{r}_j} - \underline{\delta}_{\vec{w}_j'}  \nonumber \\
		&= \underline\Delta_{\vec{w}_{j-1}'} + \alpha_{j-1} \delta_{\vec{u}_{j-1}} + A \delta_{\vec{r}_j} - \underline{\delta}_{\vec{w}_j'}\nonumber \\
		&= \underline\Delta_{\vec{w}_0'} + \sum_{\ell=1}^{j} \left[ \alpha_{\ell-1} \delta_{\vec{u}_{\ell-1}} + \vec{A}\delta_{\vec{r}_{\ell}} + \delta_{\vec{w}_\ell'} \right].
			\label{eq:Dwnr}
\end{align}

Thus, comparing~\cref{eq:Dwr} with~\cref{eq:Dwnr}, we see that if we use recomputation, \( \Delta_{\vec{w}_j'} \) is essentially ``reset'' in each iteration; the errors made do not accumulate and \( \Delta_{\vec{w}_j'} \) remains a small multiple of the machine precision.
Without recomputation, however, a bound on size of the quantity \( \underline{\Delta}_{\vec{w}_j'} \) will be monotonically increasing, growing larger in each iteration due to the accumulation of rounding errors. 
We note that the expression for the maximum attainable accuracy for \PipePRCG in \cref{eq:Dwnr}, \emph{without the recomputation of \(\vec{w}_j\)}, closely resembles the expression derived for the maximum attainable accuracy for \GVCG derived in~\cite{cools_yetkin_agullo_giraud_vanroose_18}.
This is demonstrated numerically in \cref{fig:residual_gap-three_term_error}.
Indeed, without recomputation, it is easily observed that \PipePRCG and \GVCG typically behave quite similarly.

If we look at the expression for \( \Delta_{\vec{s}_{i}} \), we see that \( \Delta_{\vec{w}_j'} \) (or in the case of no recomputation, \( \underline{\Delta}_{\vec{w}_j'} \)), along with local rounding errors \( \vec{A}\delta_{\vec{p}_j} \) and \( \delta_{\vec{s}_j} \), is amplified by a product of \( \beta \) terms, \( \prod_{\ell=j+1}^{i} \beta_\ell \).
In exact arithmetic at least, we have 
\begin{align*}
\prod_{\ell=j+1}^{i} \beta_\ell = \beta_{j+1} \beta_{j+2} \cdots \beta_{i} = \frac{\|\vec{r}_{i} \|^2}{\|\vec{r}_{j} \|^2},
\end{align*} 
and we expect this quantity to be greater than one in cases where CG exhibits large oscillation in the residual norm. 

Thus it is clear why recomputation can improve the maximal attainable accuracy: the term \( \Delta_{\vec{w}_j'} \) should remain small (some multiple of \( \epsilon \)) in this case, whereas without recomputation, the analogous term can grow larger and larger throughout the iterations.
It is also clear why the \PipePRCG method can exhibit potentially decreased accuracy versus the \HSCG method; the local rounding errors (including the term \( \Delta_{\vec{w}_j'} \)) are amplified by quantities related to the ratios of residual norms, which can be large for certain problems. 

Finally, we note that the same theory holds for \PipePRMCG as the expressions for the relevant recurrences are the same as in \PipePRCG. 


\subsubsection{Towards understanding convergence}

Following the same procedure as for \PRCG, we can write
\begin{align*}
    \vec{s}_{k}
    &=
    \vec{w}_{k}' + \beta_{k} \vec{s}_{k-1} + \delta_{\vec{s}_{k}}
    = \vec{A}\vec{r}_{k} + \beta_{k} \vec{s}_{k-1} - \Delta_{\vec{w}_k'} + \delta_{\vec{s}_{k}}. 
\end{align*}
This is almost the same as \cref{eqn:prcg:sk} except that \( \vec{A}\delta_{\vec{p}_{k}} - \beta_k \delta_{\vec{s}_{k-1}} \) has been replaced by \( -\Delta_{\vec{w}_k'} \).
Thus, through similar computations as above, we find
\begin{align*}
    \vec{A}\vec{q}_{k}
    &=
    \frac{1}{\alpha_{k-1}} \frac{\|\vec{r}_k\|}{\|\vec{r}_{k-1}\|} \vec{q}_{k+1} + \left( \frac{1}{\alpha_{k-1}} + \frac{\beta_{k-1}}{\alpha_{k-2}} \right) \vec{q}_{k} + \frac{1}{\alpha_{k-2}} \frac{\|\vec{r}_{k-1}\|}{\|\vec{r}_{k-2}\|} \vec{q}_{k-1} + \vec{F}_{k},
\end{align*}
where there error term is given by
\begin{align*}
    \vec{F}_k
    &= \frac{\Delta_{\beta_{k-1}}}{\alpha_{k-2}} \frac{\|\vec{r}_{k-2}\|^2}{\|\vec{r}_{k-1}\|^2} \vec{q}_{k-1}  - \frac{(-1)^k}{\|\vec{r}_{k-1}\|} \vec{f}_k
\end{align*}
and
\begin{align*}
   \vec{f}_k 
    &= 
    - \frac{\beta_{k-1}}{\alpha_{k-2}} \delta_{\vec{r}_{k-1}} + \Delta_{\vec{w}_k'} - \delta_{\vec{s}_{k-1}} + \delta_{\vec{r}_k}.
\end{align*}

While the three-term Lanczos recurrence generated by \PipePRCG is satisfied to within local rounding errors, the error term \( \vec{F}_k \) for \PipePRCG depends on \( \Delta_{w_k'} \) rather than \( \vec{A} \delta_{\vec{p}_k} \).
Thus, it is reasonable to expect it to be larger than the error term for \PRCG.
On the other hand, the expression for the degree to which \GVCG satisfies the three-term Lanczos recurrence depends on errors made in all previous steps \cite{greenbaum_liu_chen_19}.
This provides some intuition for why the rate of convergence of \PipePRCG is observed to be better than \GVCG.
Again, a numerical comparison is shown in \cref{fig:residual_gap-three_term_error}.

\section{Numerical experiments}
\label{sec:numerical}

\begin{figure*}[ht]\centering
    \includegraphics[width=\textwidth]{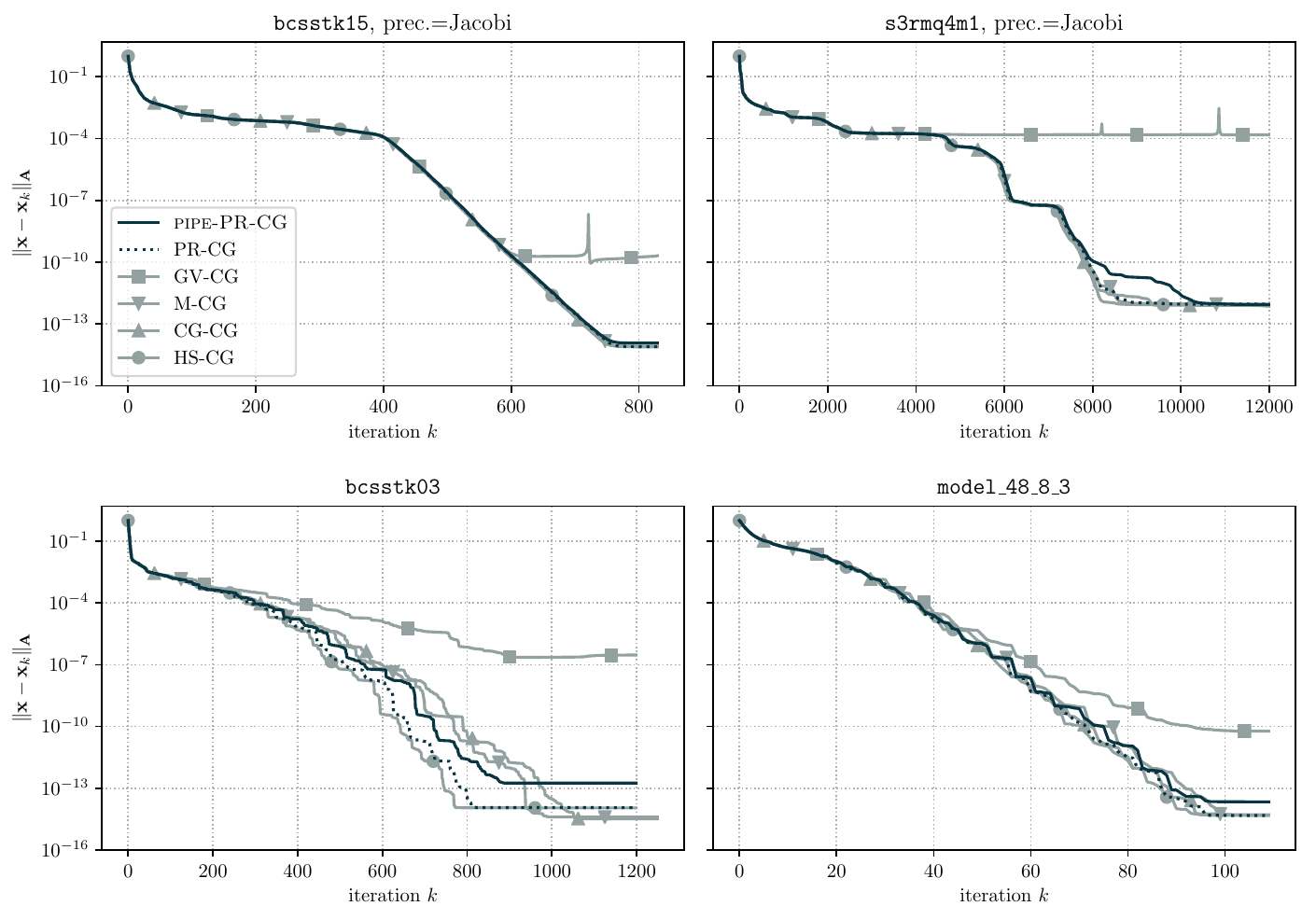}
    \caption{Error \( \vec{A} \)-norm of conjugate gradient variants on selected problems.}
    \label{fig:error_A_norm}
\end{figure*}
In this section we present the results of numerical experiments intended to give insight into the numerical behavior of the variants introduced in this paper.
We run experiments on a range of matrices from the Matrix Market \cite{boisvert_pozo_remington_barrett_dongarra_97}.
In addition we include a model problem which was introduced in \cite{strakos_91} and has since been considered in \cite{strakos_greenbaum_92,dazevedo_romine_92,greenbaum_liu_chen_19}.
This model problem has eigenvalues 
\( \lambda_1 = 1/\kappa \), \( \lambda_n = 1 \), and for \( i=2,3,\ldots,n-1 \), 
\begin{align}
    \lambda_i = \lambda_1 + \left( \frac{i-1}{n-1} \right) \cdot (\lambda_n - \lambda_1) \cdot \rho^{n-i}
		\label{eq:model}
\end{align}
for \( n=48 \), \( \rho = 0.8 \), and \( \kappa=10^3 \), with eigenvectors chosen uniformly from the set of unitary matrices.
Since the spacing between eigenvalues grows exponentially, this is a particularly difficult problem in finite precision arithmetic.

The experiments are implemented in NumPy using IEEE double precision floating point arithmetic.
The results are outlined in \cref{table:convergence}.
In this table we give two summary statistics: (i) the number of iterations required to decrease the \( \vec{A} \)-norm of the error by a factor of \( 10^{5} \), and (ii) the minimum error reached.
For a given problem, these two quantities give a rough indication of the rate of convergence and ultimately attainable accuracy.
Plots of convergence for all experiments appearing in \cref{table:convergence} can be found online in the repository linked in \cref{sec:resources}.
As was done in \cite{ghysels_vanroose_14}, the right-hand side \( \vec{b} \) is chosen so that \( \vec{x}^* = \vec{A}^{-1}\vec{b} \) has entries \( 1/\sqrt{n} \), and the initial guess \( \vec{x}_0 \) is the zero vector.
For most problems we selected, we run tests without a preconditioner, and then with a simple Jacobi (diagonal) preconditioner.
For each problem, we run all variants for a sufficient number of iterations such that the true residual stagnates.

\afterpage{\setlength{\tabcolsep}{0.5em} 
\renewcommand{\arraystretch}{1.1}
\begin{table*}\centering
\scalebox{.5}{
    \begin{tabular}{l|crr|rrrrrrr|rrrrrrr}
        \toprule
        \bmcol{l}{matrix} & 
        \bmcol{r}{prec.} & 
        \bmcol{r}{\( \boldsymbol{n} \)} & 
        \bmcol{r}{nnz} & 
        \multicolumn{7}{c}{\textbf{iterations: \( \min\{ k : \| \vec{e}_k \|_\vec{A} / \| \vec{e}_0 \|_\vec{A} < 10^{-5} \)\} }} &
        \multicolumn{7}{c}{\textbf{minimum error: \( \min_k \log_{10}( \| \vec{e}_k \|_\vec{A} / \| \vec{e}_0 \|_\vec{A} ) \)}} 
        \\ \midrule
        \multicolumn{4}{c}{} &
        \mcol{r}{\small\hyperref[alg:hs_pcg]{HS}} & \mcol{r}{\small\hyperref[alg:cg_pcg]{CG}} & \mcol{r}{\small\hyperref[alg:m_pcg]{M}} & \mcol{r}{\small\hyperref[alg:pr_pcg]{PR}} & \mcol{r}{\small\hyperref[alg:gv_pcg]{GV}} & \mcol{r}{\small\hyperref[alg:pipe_m_pcg]{\textsc{PM}}} & \mcol{r}{\small\hyperref[alg:pipe_pr_pcg]{\textsc{PPR}}} &
        \mcol{r}{\small\hyperref[alg:hs_pcg]{HS}} & \mcol{r}{\small\hyperref[alg:cg_pcg]{CG}} & \mcol{r}{\small\hyperref[alg:m_pcg]{M}} & \mcol{r}{\small\hyperref[alg:pr_pcg]{PR}} & \mcol{r}{\small\hyperref[alg:gv_pcg]{GV}} & \mcol{r}{\small\hyperref[alg:pipe_m_pcg]{\textsc{PM}}} & \mcol{r}{\small\hyperref[alg:pipe_pr_pcg]{\textsc{PPR}}}
        \\
        \cmidrule(r){5-5} \cmidrule(r){6-8} \cmidrule(r){9-11} \cmidrule(r){12-12} \cmidrule(r){13-15} \cmidrule(r){16-18}
        \texttt{1138\_bus} & - & 1138 & 4054& {1721}& {1753}& {1797}& {1727}& {1870}& {1799}& {1733}&{-12.69}&{-12.52}&{-12.70}&{-12.73}&\tableemph{-6.54}&{-11.85}&{-11.85}\\ 
\texttt{494\_bus} & - & 494 & 1666& {898}& {917}& {941}& {899}& \tableemph{1040}& {957}& {909}&{-13.14}&{-12.48}&{-13.11}&{-13.11}&\tableemph{-6.89}&{-12.24}&{-12.16}\\ 
\texttt{662\_bus} & - & 662 & 2474& {443}& {447}& {451}& {444}& {464}& {451}& {444}&{-13.93}&{-13.21}&{-13.93}&{-13.95}&\tableemph{-8.78}&{-12.67}&{-13.35}\\ 
\texttt{685\_bus} & - & 685 & 3249& {437}& {456}& {455}& {439}& \tableemph{485}& {471}& {445}&{-14.36}&{-14.10}&{-14.31}&{-14.36}&\tableemph{-9.53}&{-13.15}&{-13.06}\\ 
\texttt{bcsstk03} & - & 112 & 640& {364}& \tableemph{439}& \tableemph{425}& {380}& \tableemph{598}& \tableemph{492}& \tableemph{411}&{-14.55}&{-14.49}&{-14.40}&{-14.43}&\tableemph{-6.86}&\tableemph{-12.65}&\tableemph{-12.96}\\ 
\texttt{bcsstk14} & - & 1806 & 63454& {3982}& {4060}& {4045}& {4003}& {4212}& {4096}& {4014}&{-14.31}&{-14.19}&{-14.30}&{-14.29}&\tableemph{-5.32}&{-14.24}&{-14.26}\\ 
\texttt{bcsstk15} & - & 3948 & 117816& {5702}& {5777}& {5782}& {5721}& {5951}& {5820}& {5738}&{-13.77}&{-13.40}&{-13.79}&{-13.79}&\tableemph{-5.60}&{-13.61}&{-13.74}\\ 
\texttt{bcsstk16} & - & 4884 & 290378& {429}& {430}& {430}& {430}& \tableemph{818}& {431}& {430}&{-14.48}&{-14.34}&{-14.47}&{-14.47}&\tableemph{-5.07}&{-14.18}&{-14.20}\\ 
\texttt{bcsstk17} & - & 10974 & 428650& {17568}& {18103}& {18174}& {17590}& \tableemph{-}& {18497}& {17795}&{-13.42}&{-12.85}&{-13.42}&{-13.45}&\tableemph{-4.25}&\tableemph{-12.06}&{-12.30}\\ 
\texttt{bcsstk18} & - & 11948 & 149090& {42525}& {43366}& {43118}& {42652}& \tableemph{49210}& {43704}& {42908}&{-13.18}&{-13.15}&{-13.16}&{-13.18}&\tableemph{-5.00}&{-13.15}&{-13.18}\\ 
\texttt{bcsstk27} & - & 1224 & 56126& {519}& {521}& {523}& {520}& {531}& {523}& {520}&{-14.42}&{-14.26}&{-14.42}&{-14.43}&\tableemph{-9.85}&{-14.04}&{-14.10}\\ 
\texttt{bcsstm19} & - & 817 & 817& {274}& {287}& {286}& {274}& \tableemph{334}& {299}& {277}&{-14.84}&{-14.84}&{-14.82}&{-14.82}&\tableemph{-9.63}&{-14.59}&{-14.67}\\ 
\texttt{bcsstm20} & - & 485 & 485& {203}& {219}& {221}& {205}& \tableemph{239}& \tableemph{228}& {208}&{-15.06}&{-14.95}&{-15.09}&{-15.13}&\tableemph{-9.81}&{-14.92}&{-14.89}\\ 
\texttt{bcsstm21} & - & 3600 & 3600& {3}& {3}& {3}& {3}& {3}& {3}& {3}&{-15.69}&{-15.68}&{-15.69}&{-16.47}&{-14.58}&{-14.92}&{-14.92}\\ 
\texttt{bcsstm22} & - & 138 & 138& {43}& {43}& {43}& {43}& {43}& {43}& {43}&{-15.43}&{-15.43}&{-15.48}&{-15.43}&\tableemph{-12.51}&{-15.26}&{-15.07}\\ 
\texttt{bcsstm23} & - & 3134 & 3134& {1325}& {1376}& {1360}& {1342}& {1434}& {1367}& {1346}&{-14.35}&{-14.29}&{-14.31}&{-14.36}&\tableemph{-6.87}&{-14.30}&{-14.34}\\ 
\texttt{bcsstm24} & - & 3562 & 3562& {1573}& {1689}& {1686}& {1595}& \tableemph{19411}& {1698}& {1605}&{-14.14}&{-13.97}&{-14.05}&{-14.13}&\tableemph{-5.00}&{-13.83}&{-14.04}\\ 
\texttt{bcsstm25} & - & 15439 & 15439& {10089}& {10948}& {10963}& {10293}& \tableemph{12736}& \tableemph{11245}& {10400}&{-13.84}&{-13.74}&{-13.68}&{-13.81}&\tableemph{-5.50}&{-13.62}&{-13.76}\\ 
\texttt{model\_48\_8\_3} & - & 48 & 2304& {43}& {42}& {45}& {44}& {45}& {43}& {44}&{-14.32}&{-14.28}&{-14.31}&{-14.32}&\tableemph{-10.23}&{-13.67}&{-13.66}\\ 
\texttt{nos1} & - & 237 & 1017& {1846}& {1895}& {2008}& {1843}& \tableemph{2305}& {1999}& {1870}&{-12.81}&{-13.07}&{-12.83}&{-12.80}&\tableemph{-5.81}&{-12.15}&{-11.82}\\ 
\texttt{nos2} & - & 957 & 4137& {29829}& {30672}& {32717}& {29706}& \tableemph{-}& {32157}& {29744}&{-11.29}&{-11.21}&{-11.30}&{-11.29}&\tableemph{-3.23}&{-10.99}&{-10.99}\\ 
\texttt{nos3} & - & 960 & 15844& {221}& {221}& {221}& {221}& {221}& {221}& {221}&{-13.39}&{-13.58}&{-13.39}&{-13.39}&\tableemph{-9.86}&{-13.36}&{-13.22}\\ 
\texttt{nos4} & - & 100 & 594& {72}& {72}& {72}& {72}& {72}& {72}& {72}&{-14.33}&{-14.41}&{-14.32}&{-14.33}&\tableemph{-11.47}&{-14.19}&{-14.19}\\ 
\texttt{nos5} & - & 468 & 5172& {315}& {315}& {316}& {316}& {317}& {316}& {316}&{-14.98}&{-14.97}&{-15.00}&{-14.99}&\tableemph{-10.97}&{-14.91}&{-14.89}\\ 
\texttt{nos6} & - & 675 & 3255& {551}& {555}& {582}& {555}& \tableemph{672}& {601}& {589}&{-12.21}&{-12.28}&{-12.23}&{-12.22}&\tableemph{-6.67}&\tableemph{-10.21}&\tableemph{-10.21}\\ 
\texttt{nos7} & - & 729 & 4617& {2869}& {2798}& \tableemph{3536}& {2874}& \tableemph{-}& \tableemph{3416}& {2899}&{-9.01}&{-8.74}&{-8.97}&{-9.01}&\tableemph{-0.65}&\tableemph{-6.80}&\tableemph{-7.24}\\ 
\texttt{s1rmq4m1} & - & 5489 & 281111& {3406}& {3447}& {3432}& {3410}& {3603}& {3442}& {3434}&{-13.54}&{-13.53}&{-13.55}&{-13.56}&\tableemph{-8.23}&{-13.47}&{-13.43}\\ 
\texttt{s1rmt3m1} & - & 5489 & 219521& {3890}& {3932}& {3910}& {3895}& {4076}& {3916}& {3908}&{-13.39}&{-13.28}&{-13.39}&{-13.40}&\tableemph{-7.35}&{-13.34}&{-13.22}\\ 
\texttt{s2rmq4m1} & - & 5489 & 281111& {10476}& {10699}& {10651}& {10491}& \tableemph{11622}& {10693}& {10615}&{-13.09}&{-13.32}&{-13.11}&{-13.12}&\tableemph{-6.07}&{-12.78}&{-13.19}\\ 
\texttt{s2rmt3m1} & - & 5489 & 219521& {14484}& {14727}& {14655}& {14533}& \tableemph{-}& {14679}& {14620}&{-12.83}&{-12.70}&{-12.82}&{-12.81}&\tableemph{-4.55}&{-12.59}&{-12.78}\\ 
\texttt{s3rmq4m1} & - & 5489 & 281111& {26628}& \tableemph{29395}& {28004}& {26937}& \tableemph{-}& {28822}& {28161}&{-12.06}&{-11.94}&{-12.06}&{-12.10}&\tableemph{-4.20}&{-11.69}&{-11.59}\\ 
\texttt{s3rmt3m1} & - & 5489 & 219521& {38459}& {41037}& {40188}& {38471}& \tableemph{-}& {40839}& {40105}&{-12.04}&{-11.86}&{-12.08}&{-12.07}&\tableemph{-4.07}&{-11.54}&{-11.28}\\ 
\texttt{s3rmt3m3} & - & 5357 & 207695& {69095}& {72598}& {71471}& {69051}& \tableemph{-}& {72258}& {70852}&{-12.60}&{-11.72}&{-12.59}&{-12.71}&\tableemph{-4.20}&{-11.65}&\tableemph{-11.13}\\ 
\texttt{1138\_bus} & Jac. & 1138 & 4054& {734}& {734}& {734}& {734}& {734}& {734}& {734}&{-12.69}&{-12.75}&{-12.67}&{-12.70}&\tableemph{-8.62}&{-12.66}&{-12.65}\\ 
\texttt{494\_bus} & Jac. & 494 & 1666& {371}& {371}& {371}& {371}& {371}& {371}& {371}&{-13.15}&{-13.09}&{-13.09}&{-13.15}&\tableemph{-9.84}&{-13.14}&{-13.16}\\ 
\texttt{662\_bus} & Jac. & 662 & 2474& {166}& {166}& {166}& {166}& {166}& {166}& {166}&{-14.16}&{-14.12}&{-14.15}&{-14.19}&\tableemph{-10.94}&{-14.11}&{-13.76}\\ 
\texttt{685\_bus} & Jac. & 685 & 3249& {192}& {192}& {192}& {192}& {192}& {192}& {192}&{-14.48}&{-14.36}&{-14.59}&{-14.46}&\tableemph{-11.32}&{-14.51}&{-14.33}\\ 
\texttt{bcsstk03} & Jac. & 112 & 640& {118}& {118}& {120}& {120}& {120}& {120}& {121}&{-14.10}&{-14.11}&{-14.10}&{-14.05}&\tableemph{-9.48}&{-13.48}&{-13.50}\\ 
\texttt{bcsstk14} & Jac. & 1806 & 63454& {198}& {198}& {198}& {198}& {198}& {198}& {198}&{-14.78}&{-14.73}&{-14.66}&{-14.67}&\tableemph{-11.18}&{-14.35}&{-14.36}\\ 
\texttt{bcsstk15} & Jac. & 3948 & 117816& {442}& {442}& {443}& {444}& {444}& {444}& {444}&{-14.10}&{-14.08}&{-14.11}&{-14.10}&\tableemph{-10.05}&{-13.95}&{-13.93}\\ 
\texttt{bcsstk16} & Jac. & 4884 & 290378& {132}& {132}& {132}& {132}& {132}& {132}& {132}&{-14.61}&{-14.51}&{-14.60}&{-14.60}&\tableemph{-10.98}&{-14.24}&{-14.20}\\ 
\texttt{bcsstk17} & Jac. & 10974 & 428650& {2203}& {2205}& {2210}& {2212}& {2218}& {2214}& {2216}&{-13.98}&{-13.71}&{-13.98}&{-14.00}&\tableemph{-7.99}&{-13.40}&{-13.32}\\ 
\texttt{bcsstk18} & Jac. & 11948 & 149090& {536}& {537}& {539}& {541}& {542}& {541}& {542}&{-14.57}&{-14.54}&{-14.57}&{-14.55}&\tableemph{-10.12}&{-14.30}&{-14.30}\\ 
\texttt{bcsstk27} & Jac. & 1224 & 56126& {173}& {173}& {173}& {174}& {174}& {174}& {174}&{-14.67}&{-14.41}&{-14.67}&{-14.70}&\tableemph{-10.45}&{-13.99}&{-14.03}\\ 
\texttt{model\_48\_8\_3} & Jac. & 48 & 2304& {49}& {48}& {50}& {50}& {52}& {50}& {50}&{-14.30}&{-14.25}&{-14.28}&{-14.29}&\tableemph{-10.66}&{-13.70}&{-13.72}\\ 
\texttt{nos1} & Jac. & 237 & 1017& {306}& {314}& {322}& {312}& \tableemph{346}& {323}& {326}&{-12.98}&{-12.79}&{-12.93}&{-12.96}&\tableemph{-6.70}&{-12.67}&{-12.28}\\ 
\texttt{nos2} & Jac. & 957 & 4137& {3047}& {3184}& {3197}& {3097}& \tableemph{-}& {3326}& {3303}&{-11.27}&{-11.32}&{-11.30}&{-11.27}&\tableemph{-3.44}&{-11.23}&{-11.12}\\ 
\texttt{nos3} & Jac. & 960 & 15844& {186}& {186}& {186}& {186}& {186}& {186}& {186}&{-13.38}&{-13.42}&{-13.37}&{-13.39}&\tableemph{-9.62}&{-13.55}&{-13.53}\\ 
\texttt{nos4} & Jac. & 100 & 594& {67}& {67}& {67}& {67}& {67}& {67}& {67}&{-14.30}&{-14.39}&{-14.34}&{-14.36}&\tableemph{-11.76}&{-14.22}&{-14.14}\\ 
\texttt{nos5} & Jac. & 468 & 5172& {136}& {136}& {136}& {136}& {136}& {136}& {136}&{-15.07}&{-14.99}&{-15.11}&{-15.08}&\tableemph{-12.02}&{-14.82}&{-14.89}\\ 
\texttt{nos6} & Jac. & 675 & 3255& {71}& {71}& {71}& {71}& {71}& {71}& {71}&{-12.17}&{-12.00}&{-12.19}&{-12.20}&\tableemph{-9.10}&{-12.14}&{-12.14}\\ 
\texttt{nos7} & Jac. & 729 & 4617& {67}& {67}& {67}& {67}& {67}& {67}& {67}&{-8.91}&{-9.21}&{-8.90}&{-8.88}&\tableemph{-6.41}&{-9.42}&{-9.41}\\ 
\texttt{s1rmq4m1} & Jac. & 5489 & 281111& {595}& {595}& {596}& {596}& {597}& {597}& {597}&{-13.95}&{-13.86}&{-13.97}&{-13.90}&\tableemph{-8.49}&{-13.67}&{-13.65}\\ 
\texttt{s1rmt3m1} & Jac. & 5489 & 219521& {674}& {674}& {674}& {675}& {675}& {675}& {676}&{-13.58}&{-13.76}&{-13.60}&{-13.62}&\tableemph{-8.76}&{-13.71}&{-13.20}\\ 
\texttt{s2rmq4m1} & Jac. & 5489 & 281111& {1437}& {1437}& {1438}& {1439}& {1439}& {1439}& {1440}&{-13.16}&{-12.65}&{-13.15}&{-13.23}&\tableemph{-6.65}&{-13.10}&{-12.81}\\ 
\texttt{s2rmt3m1} & Jac. & 5489 & 219521& {2030}& {2028}& {2033}& {2034}& {2040}& {2037}& {2039}&{-12.84}&{-12.65}&{-12.83}&{-12.81}&\tableemph{-6.19}&{-12.47}&{-12.19}\\ 
\texttt{s3dkq4m2} & Jac. & 90449 & 4820891& {25527}& {25513}& {25548}& {25553}& \tableemph{-}& {25576}& {25582}&{-11.09}&{-11.17}&{-11.10}&{-11.09}&\tableemph{-4.03}&{-11.12}&{-11.21}\\ 
\texttt{s3dkt3m2} & Jac. & 90449 & 3753461& {36195}& {36152}& {36247}& {36263}& \tableemph{-}& {36327}& {36348}&{-11.39}&{-10.91}&{-11.39}&{-11.39}&\tableemph{-3.97}&{-11.09}&{-10.82}\\ 
\texttt{s3rmq4m1} & Jac. & 5489 & 281111& {5743}& {5726}& {5775}& {5780}& \tableemph{-}& {5800}& {5806}&{-12.06}&{-12.14}&{-12.03}&{-12.04}&\tableemph{-3.82}&{-12.00}&{-12.06}\\ 
\texttt{s3rmt3m1} & Jac. & 5489 & 219521& {8827}& {8806}& {8867}& {8871}& \tableemph{-}& {8908}& {8917}&{-12.07}&{-11.78}&{-12.09}&{-12.10}&\tableemph{-3.78}&{-11.57}&{-11.41}\\ 
\texttt{s3rmt3m3} & Jac. & 5357 & 207695& {10251}& {10248}& {10317}& {10324}& \tableemph{-}& {10385}& {10404}&{-12.84}&{-11.77}&{-12.58}&{-12.93}&\tableemph{-4.43}&{-11.66}&{-11.95}\\ 
 
    \bottomrule
    \end{tabular}
}
    \begin{minipage}{\textwidth}
    \vspace{1em}
        \caption{
        Summary statistics of convergence behavior on problems from the Matrix Market.
        Preconditioners are applied using preconditioned variants rather than constructing an explicitly-preconditioned system.
        Values are bold if they differ from \HSCG by more than ten percent, and dashes indicate that a method failed to reach the specified accuracy.
        Note that \( \vec{e}_k = \vec{A}^{-1} \vec{b} - \vec{x}_k \) is the error at step \( k \).
        }
        \label{table:convergence}
    \end{minipage}
\end{table*}
\clearpage}

\cref{fig:error_A_norm} shows the results of four of the numerical tests contained in \cref{table:convergence}.
These problems were chosen to highlight some of the types of behavior observed for finite precision CG variants.
On some problems such as {\tt bcsstk03} and the model problem {\tt model\_48\_8\_3}, the rate of convergence and final accuracy of each variant may differ, due primarily to the large gaps in the spectrum \cite{carson_rozloznik_strakos_tichy_tuma_18}.
Alternatively, on many other problems, such as {\tt bcsstk15} with Jacobi preconditioning, the rate of convergence for all variants is the same until the final accuracy is reached.
However, even on such problems, it may be possible for the final accuracy of a variant to be significantly worse than other variants. 
For instance, on \texttt{s3rmq4m1} with Jacobi preconditioning, the final accuracy of \GVCG is \emph{8 orders of magnitude} worse than \HSCG even though the rates of convergence are initially the same.

We note that on problems where \CGCG encounters a delay of convergence, such as \texttt{bcsstk03}, \PRCG converges more quickly.
More notably, for the experiments in \cref{table:convergence}, the pipelined predict-and-recompute variants \PipePRMCG and \PipePRCG show significantly better convergence than \GVCG, frequently exhibiting convergence similar to that of \HSCG.
In particular, on all the problems tested, \PipePRMCG and \PipePRCG converge to a final accuracy within 10 percent (on a log scale) of that of \HSCG if Jacobi preconditioning is used, and on some problems, these two variants actually converge to a \emph{better} final accuracy than \HSCG.

\begin{figure*}[t]\centering
    \begin{subfigure}{.48\textwidth}
        \includegraphics[width=\textwidth]{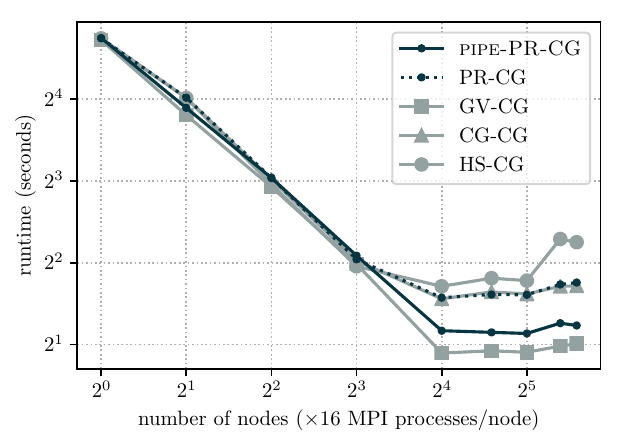}
    \end{subfigure}\hfill
    \begin{subfigure}{.48\textwidth}
        \includegraphics[width=\textwidth]{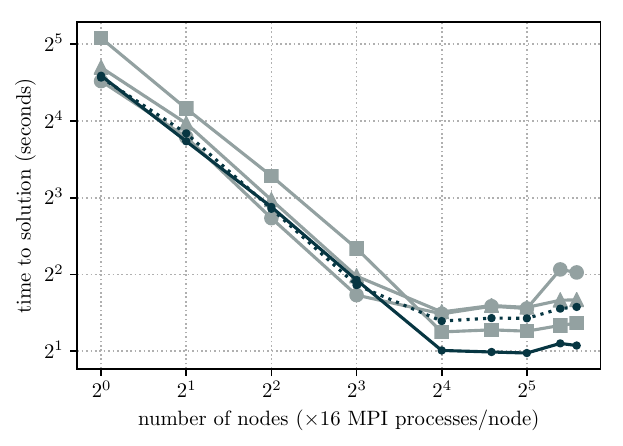}
    \end{subfigure}
    \caption{Strong scaling experiment on (dense) model problem with \( n=10752 \) unknowns and \( 1000 \) iterations.
    Time to solution is approximate time to decrease the \( \vec{A} \)-norm of error by a factor of \( 10^5 \). }
    \label{fig:strong_scale}
\end{figure*}

\section{Parallel performance}

We implement the unpreconditioned variants discussed in this paper using mpi4py \cite{mpi4py} and perform a strong scaling experiment on the Hyak supercomputer at the University of Washington.
In this experiment we solve the model problem \cref{eq:model} with parameters \( n=10752 \), \( \rho = 0.9 \), \( \kappa = 10^6 \).
In order to demonstrate that the additional computation costs associated with the extra matrix-vector product required by \PipePRCG are not particularly important, we represent our model problem as a dense matrix.
In each trial, we first allocate 48 nodes with 16 processors each and then iterate over a selection of node counts.
On each node count we run each variant with for \( 1000 \) iterations timing only the main loop of each variant and not any costs corresponding to setup.
In order to account for effects such as cluster topology and system noise we repeat this process 4 times and report the minimum runtimes and accuracy of these runs.

\Cref{fig:strong_scale} shows the results of the experiment.
In our implementation, we compute \( \vec{A} \vec{s}_k \) and \( \vec{A} \vec{r}_k \) simultaneously using a block vector \( [\vec{s}_k,\vec{r}_k] \) so that \( \vec{A} \) is accessed only once per iteration.
While the number of floating point operations is nearly doubled in \PipePRCG, the memory access pattern is more or less unchanged.
Specifically, \( \vec{A} \) needs to only be read/communicated a single time per iteration. 
For this reason, the runtime of \PipePRCG is almost the same as \GVCG.
If we compute the two matrix products separately, then the runtime of \PipePRCG is approximately doubled on low node counts.

In both cases, as a higher number of nodes are used, and communication costs begin to come into play, \PipePRCG becomes faster than the non-pipelined variants and exhibits similar scaling properties as \GVCG.
In the high node limit, \GVCG does perform marginally better than \PipePRCG. 
This is because our test involves a dense matrix, so due to the two matrix-vector products, the amount of data communicated between nodes in \PipePRCG is nearly double that of \GVCG.
However, it is important to keep in mind that the actual ``time to solution'' is given by the product of the time per iteration and the number of iterations required.

\section{Future work}
\label{sec:future}

First, while we have provided an analysis of the maximal attainable accuracy for our new variants, the analysis of the resulting convergence rate is far from complete. This is a difficult area of study and deserves further rigorous treatment. 

Second, since the maximum reduction in communication cost of \PipePRCG over \HSCG is only a factor of three, potential ways of further decreasing communication costs should be explored.
Recently, there has been work on ``deep pipelined'' CG variants where more matrix-vector products are overlapped with global communication; see, for instance, \cite{ghysels_vanroose_14,cornelis_cools_vanroose_19,cools_cornelis_vanroose_19}.
This approach is similar to the ``look ahead'' strategy suggested in \cite{rosendale_83} and may be compatible with the schemes introduced here.

Alternatively, it may be possible to either incorporate predict-and-recompute strategies into \( s \)-step methods or to develop new \( s \)-step methods which are built on \PRCG.
As \CGCG is the \( s=1 \) case of the \( s \)-step method from \cite{chronopoulos_gear_89}, we expect it may be possible to develop an \( s \)-step method based based on \PRCG which has slightly better numerical properties than \CGCG.
Finding such a method which is usable in practice would be of great practical interest. 



Finally, the predict-and-recompute variants presented here can be naturally extended to other related methods such as conjugate residual and conjugate gradient squared.

\section{Conclusion}
In this paper we extend the predict-and-recompute idea of Meurant \cite{meurant_87} and derive new pipelined CG variants, \PipePRMCG and \PipePRCG. 
These variants exhibit better theoretical scaling properties than the standard \HSCG algorithm as well as improved numerical behavior compared to their previously-studied counterpart, \GVCG, on \emph{every} numerical experiment we ran. 
We provide an analysis of the maximum attainable accuracy, applicable to both \PipePRMCG and \PipePRCG, which explains the benefit of the predict-and-recompute with regards to avoiding the breakdown observed in \cite{meurant_87}, as well as to the improved final accuracy and rate of convergence observed by our pipelined variants.
We additionally provide a strong scaling experiment which confirms the potential performance benefits of our approach. 
Our predict-and-recompute variants require exactly the same input parameters as \HSCG, and therefore have the potential to be used wherever \HSCG is used without any additional parameter selection.
Despite these advances, there is still significant room for future work on high performance CG variants, especially in the direction of further decreasing the communication costs.

\section*{Acknowledgments}
The first author expresses gratitude to their advisor, Anne Greenbaum, for encouraging them to pursue the side projects which gave rise to the results of this paper.
We also thank everyone who provided feedback on early drafts, especially Siegfried Cools, Pieter Ghysels, and G\'erard Meurant, as well as the referees whose feedback significantly improved several aspects of this paper.

\appendix
\section{Additional resources}
\label{sec:resources}

A repository with the code necessary to reproduce all of the figures and results in this paper is available at \href{https://github.com/tchen01/new_cg_variants}{\texttt{https://github.com/tchen01/new\_cg\_variants}}, and released to the public domain under the \href{https://opensource.org/licenses/MIT}{MIT License}.
The repository also contains convergence data and plots for all the matrices listed in \cref{table:convergence}.

We are committed to facilitating the reproducibility process, and encourage questions and inquiries into the methods used in this paper.

\section{Bound on sum of products}

Define
\begin{align*}
    C = \sum_{i=1}^{n} \prod_{j=1}^{b_i} c_{i,j}.
\end{align*}

Then
\begin{align*}
    \left| \operatorname{fp} \left( \prod_{j=1}^{b_i} c_{i,j} \right)
    - \prod_{j=1}^{b_i} c_{i,j} \right|
    \leq (1+\epsilon)^{b_i} \prod_{j=1}^{b_i} \left| c_{i,j} \right|,
\end{align*}
so
\begin{align}
    |\operatorname{fp}(C) - C| 
    \leq n\, \epsilon \sum_{i=1}^{n} (1+\epsilon)^{b_i} \prod_{j=1}^{b_i} |c_{i,j}|
    =  n\, \epsilon \sum_{i=1}^{n} \prod_{j=1}^{b_i} |c_{i,j}| + \mathcal{O}(\epsilon^2).
    \label{eqn:sum_prod}
\end{align}

\section{A note on other variants studied}
\label{sec:other_variants}

Many of the scalar quantities appearing in the variants discussed in this paper can be replaced with equivalent expressions. 
For instance, the difference between \MCG and \PRCG (as well as \PipePRMCG and \PipePRCG) is the choice of expression for \( \nu_k \), which in exact arithmetic is equal to \( \langle \tilde{\vec{r}}_k, \vec{r}_k \rangle \).
However, it is not typically clear which scalar expression should be used.

We first discuss our simplified expression for \cref{eqn:nuk} which is used in \PRCG and \PipePRCG.
Recall that in exact arithmetic, \( \langle \tilde{\vec{s}}_k, \vec{r}_k \rangle = \langle \vec{s}_k, \tilde{\vec{r}}_k \rangle
\).
This cannot be expected to be true in finite precision. 

Writing the finite precision recurrences
\begin{align*}
    \tilde{\vec{s}}_k &= \vec{M}^{-1} \vec{s}_{k} + \delta_{\tilde{\vec{s}}_k}
    ,&&&
    \tilde{\vec{r}}_k &= \tilde{\vec{r}}_{k-1} - \alpha_{k-1} \tilde{\vec{s}}_{k-1} + \delta_{\tilde{\vec{r}}_k},
\end{align*}
we have
\begin{align*}
    \langle \tilde{\vec{s}}_k , \vec{r}_k \rangle
    &= \langle \vec{M}^{-1} \vec{s}_k + \delta_{\tilde{\vec{s}}_k}, \vec{r}_{k} \rangle
    \\&= \langle \vec{s}_k , \vec{M}^{-1} \vec{r}_k \rangle + \langle \delta_{\tilde{\vec{s}}_k}, \vec{r}_{k} \rangle
    \\&= \langle \vec{s}_k, \tilde{\vec{r}}_k \rangle + \langle \vec{s}_k, \Delta_{\tilde{\vec{r}}_k} \rangle + \langle \delta_{\tilde{\vec{s}}_k}, \vec{r}_k \rangle,
\end{align*}
where we have defined the \( \tilde{\vec{r}} \)-gap as \( \Delta_{\tilde{\vec{r}}_k} := \vec{M}^{-1} \vec{r}_k - \tilde{\vec{r}}_k \).
For the \( \tilde{\vec{r}} \)-gap, we have 
\begin{align*}
    \Delta_{\tilde{\vec{r}}_k} 
    &= \vec{M}^{-1} ( \vec{r}_{k-1} - \alpha_{k-1} \vec{s}_{k-1} + \delta_{\vec{r}_k} ) - ( \tilde{\vec{r}}_{k-1} - \alpha_{k-1} \tilde{\vec{s}}_{k-1} + \delta_{\tilde{\vec{r}}_k} )
    \\&= (\vec{M}^{-1} \vec{r}_{k-1} - \tilde{\vec{r}}_{k-1}) - \alpha_{k-1} (\vec{M}^{-1} \vec{s}_{k-1} - \tilde{\vec{s}}_{k-1} ) + \vec{M}^{-1} \delta_{\vec{r}_k} - \delta_{\tilde{\vec{r}}_k}
    \\&= \Delta_{\tilde{\vec{r}}_{k-1}} + \alpha_{k-1} \delta_{\tilde{\vec{s}}_{k-1}} + \vec{M}^{-1} \delta_{\vec{r}_{k}} - \delta_{\tilde{\vec{r}}_{k}}.
    \\&= \Delta_{\tilde{\vec{r}}_0} + \sum_{i=1}^{k} \left[ \alpha_{i-1} \delta_{\tilde{\vec{s}}_{i-1}} + \vec{M}^{-1} \delta_{\vec{r}_{i}} - \delta_{\tilde{\vec{r}}_{i}} \right].
\end{align*}

We note that for preconditioners with a simple structure, the errors induced by the preconditioner application may not be that large.
In such cases it is reasonable to use \( \langle \tilde{\vec{r}_k} ,  \vec{s}_k \rangle \) and \( \langle \vec{r}_k, \tilde{\vec{s}}_k \rangle \) interchangeably.
On the other hand, if the structure of the preconditioner is such that the size of the rounding errors associated with its application may be large, then it may advantageous to use the original expression \cref{eqn:nuk} rather than the reduced expression \cref{eqn:nuk_simplified}.

We now discuss the simplification of \cref{eqn:nuk} used in \MCG and \PipePRMCG.
This expression relies on the exact arithmetic relation \( \langle \vec{s}_k, \tilde{\vec{r}}_k \rangle = \langle \vec{p}_k, \vec{s}_k \rangle \).
Using the finite precision recurrence \( \alpha_k = \nu_k / \mu_k + \delta_{\alpha_k} \) we compute
\begin{align*}
    \alpha_{k} \langle \tilde{\vec{r}}_{k} , \vec{s}_{k} \rangle
    &= \left( \frac{\nu_k}{\mu_k} + \delta_{\alpha_k} \right) \langle \tilde{\vec{r}}_k, \vec{s}_k \rangle
    = \left( \frac{\langle \tilde{\vec{r}}_{k}, \vec{r}_{k} \rangle + \delta_{\nu_k}}{\langle \vec{p}_k, \vec{s}_k \rangle + \delta_{\mu_k}} + \delta_{\alpha_k} \right)  \langle \tilde{\vec{r}}_k, \vec{s}_k \rangle.
\end{align*}

In exact arithmetic, \( \langle \vec{p}_k \vec{s}_k \rangle = \langle \tilde{\vec{r}}_k, \vec{s}_k \rangle \). 
However, in finite precision
\begin{align*}
    \langle \vec{p}_k, \vec{s}_k \rangle
    &= \langle \tilde{\vec{r}}_k + \beta_{k} \vec{p}_{k-1} + \delta_{\vec{p}_k} , \vec{s}_k \rangle
    \\&= \langle \tilde{\vec{r}}_k, \vec{s}_k \rangle + \beta_k \langle \vec{p}_{k-1} + \delta_{\vec{p}_k} , \vec{A} \vec{p}_{k} + \delta_{\vec{s}_k} \rangle
    \\&= \langle \tilde{\vec{r}}_k, \vec{s}_k \rangle + \beta_k \langle \vec{p}_{k-1}, \vec{A} \vec{p}_{k} \rangle +  \langle \delta_{\vec{p}_k},  \vec{A} \vec{p}_k +  \delta_{\vec{s}_k} \rangle + \langle \delta_{\vec{p}_k}, \delta_{\vec{s}_k} \rangle .
\end{align*}

This differs from the expression for \( \langle \tilde{\vec{s}}_k, \vec{r}_k \rangle \) in that it relies on understanding the size of \( \langle \vec{p}_{k-1}, \vec{A} \vec{p}_k \rangle \) which should be small due to the orthogonality of certain vectors.
While expressions for local orthogonality of such vectors exist for \HSCG \cite[Proposition 5.19]{meurant_06}, these vectors need not be particularly orthogonal in other variants.
Indeed, it is easy to verify that in these variants such local orthogonality can be quite weak.

Finally, we make two remarks about previously-studied communication-hiding variants.
First, the formula for \( \mu_k \) in \CGCG and \GVCG is derived in a similar way to \cref{eqn:nuk} but then simplified. 
If left unsimplified, convergence of \CGCG appears to be slightly improved on difficult problems. 
However, using the unsimplified expression in \GVCG does not seem to lead to improved convergence.

Second, the predict-and-recompute idea may be applied to \GVCG by recomputing \( \vec{s}_{k} = \vec{A} \vec{p}_{k} \) and/or \( \vec{u_k} = \vec{A} \vec{s}_k \) or any of the recursively computed inner products.
While recomputing both \( \vec{s}_k \) and \( \vec{u}_k \) seems to result in even worse final accuracy, recomputing one or the other leads to a better rate of convergence and final accuracy than \GVCG on some problems (although neither is as good as \PipePRCG).
Recomputing inner products was not observed to have a significant effect.
Further study may be of value.

As an afterthought, we suggest that it may be possible to procedurally generate mathematically equivalent CG variants, and then automatically check if they have improved convergence properties. 
Perhaps, by finding many variants which work well, the similarities between them could provide insights into necessary properties for a good finite precision CG variant.

\section{Previously studied communication-hiding variants}
\label{sec:previous}

While we omitted the full descriptions of \MCG \cite{meurant_87}, \PipePRMCG, \CGCG \cite{chronopoulos_gear_89}, and \GVCG \cite{ghysels_vanroose_14} in the main paper, we include them here for completeness.
We additionally include the initialization required to fully implement these variants.
Note that not all variants require all of the variables defined in \cref{alg:initialize}.
For instance, \HSCG does not use any of the variables after \( \alpha_0 \).

\begin{algorithm}[H]
\caption{Initilizations}\label{alg:initialize}
\fontsize{10}{10}\selectfont
\begin{algorithmic}[1]
\Procedure{initialize}{}
\State{%
    \( \vec{r}_0 = \vec{b}-\vec{A}\vec{x}_0 \), 
    \( \tilde{\vec{r}}_0 = \vec{M}^{-1} \vec{r}_0 \), 
    \( \nu_0 = \langle \vec{r}_0,\tilde{\vec{r}}_0 \rangle \), 
    \( \vec{p}_0 = \tilde{\vec{r}}_0 \), 
    \( \vec{s}_0 = \vec{A} \vec{p}_0 \), 
    \( \alpha_0 = \nu_0 / \langle \vec{p}_0,\vec{s}_0 \rangle \)
    \( \vec{w}_0 = \vec{A} \tilde{\vec{r}}_0 \),
    \( \tilde{\vec{w}}_0 = \vec{M}^{-1} \vec{w}_0 \),
    \( \tilde{\vec{s}}_0 = \vec{M}^{-1} \vec{s}_0 \),
    \( \vec{u}_0 = \vec{A}\tilde{\vec{s}}_0 \),
    \( \tilde{\vec{u}}_0 = \vec{M}^{-1} \vec{u}_0 \),
    \( \tilde{\vec{s}}_0 = \vec{M}^{-1}s_0 \),
    \( \alpha_0 = \nu_0 / \langle \vec{p}_0, \vec{s}_0 \rangle \),
    \( \sigma_0 = \langle \tilde{\vec{r}}_0, \vec{s}_0 \rangle \),
    \( \gamma_0 = \langle \tilde{\vec{s}}_0, \vec{s}_0 \rangle \)
}
\EndProcedure
\end{algorithmic}
\end{algorithm}

\begin{landscape}
\begin{center}
\scalebox{.93}{
    \begin{minipage}[t]{.8\textwidth}
    \begin{algorithm}[H]
\caption{Meurant conjugate gradient}\label{alg:m_pcg}
\fontsize{10}{10}\selectfont
\begin{algorithmic}[1]
\Procedure{\MCG}{$\vec{A}$, $\vec{M}$, $\vec{b}$, $\vec{x}_0$}
\State \INIT
\For {\( k=1,2,\ldots \)}
    \State \( \vec{x}_k = \vec{x}_{k-1} + \alpha_{k-1} \vec{p}_{k-1} \)
    \State \( \vec{r}_k = \vec{r}_{k-1} - \alpha_{k-1} \vec{s}_{k-1} \)%
    ,~ \( \tilde{\vec{r}}_k = \tilde{\vec{r}}_{k-1} - \alpha_{k-1} \tilde{\vec{s}}_{k-1} \)
    \State \( \nu'_{k} = -\nu_{k-1} + \alpha_{k-1}^2 \gamma_{k-1}  \)
    \State \( \beta_{k} = \nu'_{k} / \nu_{k-1} \)
    \State \( \vec{p}_k = \tilde{\vec{r}}_k + \beta_k \vec{p}_{k-1} \)
    \State \( \vec{s}_k = \vec{A} \vec{p}_{k} \)%
    ,~ \( \tilde{\vec{s}}_k = \vec{M}^{-1} \vec{s}_k \)
    \State \( \mu_k = \langle \vec{p}_k,\vec{s}_k \rangle \)%
    ,~ \( \gamma_k = \langle \tilde{\vec{s}}_k, \vec{s}_k \rangle \)%
    ,~ \( \nu_k = \langle \tilde{\vec{r}}_k, \vec{r}_k \rangle \) 
    \State \( \alpha_k = \nu_k / \mu_k \)
\EndFor
\EndProcedure

\end{algorithmic}
\end{algorithm}

\end{minipage}
}\hfill
\scalebox{.93}{
    \begin{minipage}[t]{.8\textwidth}
    \begin{algorithm}[H]
\caption{Chronopoulos and Gear conjugate gradient}\label{alg:cg_pcg}
\fontsize{10}{10}\selectfont
\begin{algorithmic}[1]
\Procedure{\CGCG}{$\vec{A}$, $\vec{M}$, $\vec{b}$, $\vec{x}_0$}
\State \INIT
\For {\( k=1,2,\ldots \)}
    \State \( \vec{x}_k = \vec{x}_{k-1} + \alpha_{k-1} \vec{p}_{k-1} \)
    \State \( \vec{r}_k = \vec{r}_{k-1} - \alpha_{k-1} \vec{s}_{k-1} \)%
    ,~ \( \tilde{\vec{r}}_k = \vec{M}^{-1} \vec{r}_k \)
    \State \( \vec{w}_k = \vec{A} \tilde{\vec{r}}_k \)
    \State \( \nu_{k} = \langle \tilde{\vec{r}}_k,\vec{r}_k \rangle \)%
    ,~ \( \eta_k = \langle \tilde{\vec{r}}_k, \vec{w}_k\rangle \)
    \State \( \beta_k = \nu_k / \nu_{k-1} \)
    \State \( \vec{p}_k = \tilde{\vec{r}}_k + \beta_k \vec{p}_{k-1} \)
    \State \( \vec{s}_k = \vec{w}_k + \beta_k \vec{s}_{k-1} \)
    \State \( \mu_k = \eta_k - (\beta_k/\alpha_{k-1}) \nu_k \)
    \State \( a_k = \nu_k / \mu_k \) 
\EndFor
\EndProcedure

\end{algorithmic}
\end{algorithm}

\end{minipage}
}
\end{center}
\vfill
\begin{center}
\scalebox{.93}{
    \begin{minipage}[t]{.8\textwidth}
    \begin{algorithm}[H]
\caption{Pipelined predict-and-recompute Meurant conjugate gradient}\label{alg:pipe_m_pcg}
\fontsize{10}{10}\selectfont
\begin{algorithmic}[1]
\Procedure{\PipePRMCG}{$\vec{A}$, $\vec{M}$, $\vec{b}$, $\vec{x}_0$}
\State \INIT
\For {\( k=1,2,\ldots \)}
    \State \( \vec{x}_k = \vec{x}_{k-1} + \alpha_{k-1} \vec{p}_{k-1} \)
    \State \( \vec{r}_k = \vec{r}_{k-1} - \alpha_{k-1} \vec{s}_{k-1} \)
    ,~  \( \tilde{\vec{r}}_k = \tilde{\vec{r}}_{k-1} - \alpha_{k-1} \tilde{\vec{s}}_{k-1} \)
    \State \( \vec{w}'_k = \vec{w}_{k-1} - \alpha_{k-1} \vec{u}_{k-1} \)%
    ,~  \( \tilde{\vec{w}}'_k = \tilde{\vec{w}}_{k-1} - \alpha_{k-1} \tilde{\vec{u}}_{k-1} \)
    \State \( \nu'_{k} = -\nu_{k-1} + \alpha_{k-1}^2 \gamma_{k-1}  \) 
    \State \( \beta_{k} = \nu'_{k} / \nu_{k-1} \)
    \State \( \vec{p}_k = \tilde{\vec{r}}_k + \beta_k \vec{p}_{k-1} \)
    \State \( \vec{s}_k = \vec{w}'_k + \beta_k \vec{s}_{k-1} \)%
    ,~  \( \tilde{\vec{s}}_k = \tilde{\vec{w}}'_k + \beta_k \tilde{\vec{s}}_{k-1} \)
    \State \( \vec{u}_k = \vec{A} \tilde{\vec{s}}_k \)%
    ,~  \( \tilde{\vec{u}}_k = \vec{M}^{-1} \vec{u}_k \)
    \State \( \vec{w}_k = \vec{A} \tilde{\vec{r}}_k  \)
    ,~ \( \tilde{\vec{w}}_k = \vec{M}^{-1} \vec{w}_k \) 
    \State \( \mu_k = \langle \vec{p}_k, \vec{s}_k \rangle \)
    ,~ \( \gamma_k = \langle \tilde{\vec{s}}_k, \vec{s}_k \rangle \)%
    ,~ \( \nu_k = \langle \tilde{\vec{r}}_k, \vec{r}_k\rangle \)
    \State \( \alpha_k = \nu_k / \mu_k \)
\EndFor
\EndProcedure
\end{algorithmic}
\end{algorithm}

\end{minipage}
}\hfill
\scalebox{.93}{
    \begin{minipage}[t]{.8\textwidth}
    \begin{algorithm}[H]
    \caption{Ghysels and Vanroose conjugate gradient}\label{alg:gv_pcg}
\fontsize{10}{10}\selectfont
\begin{algorithmic}[1]
\Procedure{\GVCG}{$\vec{A}$, $\vec{M}$, $\vec{b}$, $\vec{x}_0$}
\State \INIT
\For {\(k=1,2,\ldots\)}
    \State \( \vec{x}_k = \vec{x}_{k-1} + \alpha_{k-1} \vec{p}_{k-1} \)
    \State \( \vec{r}_k = \vec{r}_{k-1} - \alpha_{k-1} \vec{s}_{k-1} \)%
    ,~ \( \tilde{\vec{r}}_k = \tilde{\vec{r}}_{k-1} - a_{k-1} \tilde{\vec{s}}_{k-1} \)
    \State \( \vec{w}_k = \vec{w}_{k-1} - a_{k-1} \vec{u}_{k-1} \)%
    ,~ \( \tilde{\vec{w}}_k = \vec{M}^{-1} \vec{w}_k \)
    \State \( \nu_{k} = \langle \tilde{\vec{r}}_k,\vec{r}_k \rangle \)%
    ,~ \( \eta_k = \langle \tilde{\vec{r}}_k, \vec{w}_k\rangle \)
    \State \( \vec{t}_k = \vec{A}\tilde{\vec{w}}_k \)
    \State \( \beta_k = \nu_k / \nu_{k-1} \)
    \State \( \vec{p}_k = \tilde{\vec{r}}_k + \beta_k \vec{p}_{k-1} \)
    \State \( \vec{s}_k = \vec{w}_k + \beta_k \vec{s}_{k-1} \)%
    ,~ \( \tilde{\vec{s}}_k = \tilde{\vec{w}}_k + \beta_k \tilde{\vec{s}}_{k-1} \)
    \State \( \vec{u}_k = \vec{t}_k + \beta_k \vec{u}_{k-1} \)
    \State \( \mu_k = \eta_k - ( \beta_k / \alpha_{k-1} ) \nu_k \)
    \State \( \alpha_k = \nu_k / \mu_k \) 
\EndFor
\EndProcedure

\end{algorithmic}
\end{algorithm}

\end{minipage}
}
\end{center}
\end{landscape}


\bibliographystyle{siamplain}
\bibliography{predict_and_recompute}

\begin{thebibliography}{10}

\bibitem{ashby_ghysels_heirman_vanroose_12}
{\sc T.~J. Ashby, P.~Ghysels, W.~Heirman, and W.~Vanroose}, {\em The impact of
  global communication latency at extreme scales on krylov methods}, in
  Algorithms and Architectures for Parallel Processing, Y.~Xiang,
  I.~Stojmenovic, B.~O. Apduhan, G.~Wang, K.~Nakano, and A.~Zomaya, eds.,
  Berlin, Heidelberg, 2012, Springer Berlin Heidelberg, pp.~428--442.

\bibitem{ballard_carson_demmel_hoemmen_knight_schwartz_14}
{\sc G.~Ballard, E.~C. Carson, J.~Demmel, M.~Hoemmen, N.~Knight, and
  O.~Schwartz}, {\em Communication lower bounds and optimal algorithms for
  numerical linear algebra}, Acta Numerica, 23 (2014), p.~1–155.

\bibitem{boisvert_pozo_remington_barrett_dongarra_97}
{\sc R.~F. Boisvert, R.~Pozo, K.~Remington, R.~F. Barrett, and J.~J. Dongarra},
  {\em Matrix market: a web resource for test matrix collections}, in Quality
  of Numerical Software, Springer, 1997, pp.~125--137.

\bibitem{carson_rozloznik_strakos_tichy_tuma_18}
{\sc E.~C. Carson, M.~Rozložník, Z.~Strakoš, P.~Tichý, and M.~Tůma}, {\em
  The numerical stability analysis of pipelined conjugate gradient methods:
  Historical context and methodology}, SIAM Journal on Scientific Computing, 40
  (2018), pp.~A3549--A3580.

\bibitem{chronopoulos_87}
{\sc A.~Chronopoulos}, {\em A Class of Parallel Iterative Methods Implemented
  on Multiprocessors}, PhD thesis, Chicago, IL, USA, 1987.

\bibitem{chronopoulos_gear_89}
{\sc A.~Chronopoulos and C.~W. Gear}, {\em $s$-step iterative methods for
  symmetric linear systems}, Journal of Computational and Applied Mathematics,
  25 (1989), pp.~153 -- 168.

\bibitem{cools_cornelis_vanroose_19}
{\sc S.~Cools, J.~Cornelis, and W.~Vanroose}, {\em Numerically stable
  recurrence relations for the communication hiding pipelined conjugate
  gradient method}, IEEE Transactions on Parallel and Distributed Systems, 30
  (2019), pp.~2507--2522.

\bibitem{cools_yetkin_agullo_giraud_vanroose_18}
{\sc S.~Cools, E.~Fatih~Yetkin, E.~Agullo, L.~Giraud, and W.~Vanroose}, {\em
  Analyzing the effect of local rounding error propagation on the maximal
  attainable accuracy of the pipelined conjugate gradient method}, SIAM Journal
  on Matrix Analysis and Applications, 39 (2018), pp.~426--450.

\bibitem{cools_vanroose_17}
{\sc S.~Cools and W.~Vanroose}, {\em Numerically stable variants of the
  communication-hiding pipelined conjugate gradients algorithm for the parallel
  solution of large scale symmetric linear systems}, 2017.

\bibitem{cornelis_cools_vanroose_19}
{\sc J.~Cornelis, S.~Cools, and W.~Vanroose}, {\em The communication-hiding
  conjugate gradient method with deep pipelines}, 2019.

\bibitem{mpi4py}
{\sc L.~D. Dalcin, R.~R. Paz, P.~A. Kler, and A.~Cosimo}, {\em Parallel
  distributed computing using python}, Advances in Water Resources, 34 (2011),
  pp.~1124 -- 1139.
\newblock New Computational Methods and Software Tools.

\bibitem{dazevedo_romine_92}
{\sc E.~F. D`Azevedo and C.~H. Romine}, {\em Reducing communication costs in
  the conjugate gradient algorithm on distributed memory multiprocessors},
  (1992).

\bibitem{demmel_hoemmen_mohiyuddin_yelick_07}
{\sc J.~Demmel, M.~F. Hoemmen, M.~Mohiyuddin, and K.~A. Yelick}, {\em Avoiding
  communication in computing krylov subspaces}, Tech. Report UCB/EECS-2007-123,
  EECS Department, University of California, Berkeley, Oct 2007.

\bibitem{dongarra_heroux_luszczek_16}
{\sc J.~Dongarra, M.~A. Heroux, and P.~Luszczek}, {\em High-performance
  conjugate-gradient benchmark: A new metric for ranking high-performance
  computing systems}, The International Journal of High Performance Computing
  Applications, 30 (2016), pp.~3--10.

\bibitem{eller_gropp_16}
{\sc P.~R. Eller and W.~Gropp}, {\em Scalable non-blocking preconditioned
  conjugate gradient methods}, in Proceedings of the International Conference
  for High Performance Computing, Networking, Storage and Analysis, SC '16,
  Piscataway, NJ, USA, 2016, IEEE Press, pp.~18:1--18:12.

\bibitem{ghysels_vanroose_14}
{\sc P.~Ghysels and W.~Vanroose}, {\em Hiding global synchronization latency in
  the preconditioned conjugate gradient algorithm}, Parallel Computing, 40
  (2014), pp.~224 -- 238.

\bibitem{greenbaum_89}
{\sc A.~Greenbaum}, {\em Behavior of slightly perturbed lanczos and
  conjugate-gradient recurrences}, Linear Algebra and its Applications, 113
  (1989), pp.~7 -- 63.

\bibitem{greenbaum_97a}
{\sc A.~Greenbaum}, {\em Estimating the attainable accuracy of recursively
  computed residual methods}, SIAM Journal on Matrix Analysis and Applications,
  18 (1997), pp.~535--551.

\bibitem{greenbaum_97}
{\sc A.~Greenbaum}, {\em Iterative Methods for Solving Linear Systems}, Society
  for Industrial and Applied Mathematics, Philadelphia, PA, USA, 1997.

\bibitem{greenbaum_liu_chen_19}
{\sc A.~Greenbaum, H.~Liu, and T.~Chen}, {\em On the convergence rate of
  variants of the conjugate gradient algorithm in finite precision arithmetic},
  2019.

\bibitem{gutknecht_strakos_00}
{\sc M.~Gutknecht and Z.~Strakos}, {\em Accuracy of two three-term and three
  two-term recurrences for krylov space solvers}, SIAM Journal on Matrix
  Analysis and Applications, 22 (2000), pp.~213--229.

\bibitem{hestenes_stiefel_52}
{\sc M.~R. Hestenes and E.~Stiefel}, {\em Methods of conjugate gradients for
  solving linear systems}, vol.~49, NBS Washington, DC, 1952.

\bibitem{johnsson_84}
{\sc L.~Johnsson}, {\em Highly concurent algorithms for solving linear systems
  of equations}, in Elliptic Problem Solvers, Academic Press, 1984, pp.~105 --
  126.

\bibitem{mcmanus_johnson_cross_99}
{\sc K.~McManus, S.~Johnson, and M.~Cross}, {\em Communication latency hiding
  in a parallel conjugate gradient method}, in Eleventh International
  Conference on Domain Decomposition Methods (London, 1998), DDM. org,
  Augsburg, Citeseer, 1999, pp.~306--313.

\bibitem{meurant_87}
{\sc G.~Meurant}, {\em Multitasking the conjugate gradient method on the cray
  x-mp/48}, Parallel Computing, 5 (1987), pp.~267 -- 280.

\bibitem{meurant_06}
{\sc G.~Meurant}, {\em The Lanczos and Conjugate Gradient Algorithms}, Society
  for Industrial and Applied Mathematics, 2006.

\bibitem{meurant_strakos_06}
{\sc G.~Meurant and Z.~Strako{\v{s}}}, {\em The lanczos and conjugate gradient
  algorithms in finite precision arithmetic}, Acta Numerica, 15 (2006),
  pp.~471--542.

\bibitem{paige_71}
{\sc C.~C. Paige}, {\em The computation of eigenvalues and eigenvectors of very
  large sparse matrices.}, PhD thesis, University of London, 1971.

\bibitem{paige_80}
{\sc C.~C. Paige}, {\em Accuracy and effectiveness of the lanczos algorithm for
  the symmetric eigenproblem}, Linear Algebra and its Applications, 34 (1980),
  pp.~235 -- 258.

\bibitem{rosendale_83}
{\sc J.~V. Rosendale}, {\em Minimizing inner product data dependencies in
  conjugate gradient iteration.}, in ICPP, 1983.

\bibitem{saad_85}
{\sc Y.~Saad}, {\em Practical use of polynomial preconditionings for the
  conjugate gradient method}, SIAM Journal on Scientific and Statistical
  Computing, 6 (1985), pp.~865--881.

\bibitem{saad_89}
{\sc Y.~Saad}, {\em Krylov subspace methods on supercomputers}, SIAM Journal on
  Scientific and Statistical Computing, 10 (1989), pp.~1200--1232.

\bibitem{sleijpen_vandervorst_fokkema_94}
{\sc G.~L.~G. Sleijpen, H.~A. van~der Vorst, and D.~R. Fokkema}, {\em
  Bicgstab(l) and other hybrid bi-cg methods}, Numerical Algorithms, 7 (1994),
  pp.~75--109.

\bibitem{strakos_91}
{\sc Z.~Strakoš}, {\em On the real convergence rate of the conjugate gradient
  method}, Linear Algebra and its Applications, 154-156 (1991), pp.~535 -- 549.

\bibitem{strakos_greenbaum_92}
{\sc Z.~Strakoš and A.~Greenbaum}, {\em Open questions in the convergence
  analysis of the lanczos process for the real symmetric eigenvalue problem},
  University of Minnesota, 1992.

\bibitem{strzodka_goddeke_06}
{\sc R.~Strzodka and D.~Goddeke}, {\em Pipelined mixed precision algorithms on
  fpgas for fast and accurate pde solvers from low precision components}, in
  2006 14th Annual IEEE Symposium on Field-Programmable Custom Computing
  Machines, April 2006, pp.~259--270.

\end{thebibliography}

\end{document}